\documentclass[11pt,a4paper]{article}
\usepackage[T1]{fontenc}
\usepackage{lmodern}
\usepackage{amsmath,amssymb,amsthm,mathtools}
\usepackage[margin=2.55cm]{geometry}
\usepackage{enumitem,microtype}
\usepackage[hidelinks]{hyperref}

\usepackage{stmaryrd}
\usepackage{graphicx}
\usepackage{amscd}
\usepackage{amsbsy}
\usepackage{epsfig}
\usepackage{tikz}
\usepackage{cite}
\usepackage{indentfirst}
\usepackage{multicol}
\usepackage{wasysym}
\usepackage{bbm}
\usepackage{longtable,booktabs}
\numberwithin{equation}{section}
\newtheorem{theorem}{Theorem}[section]
\newtheorem{lemma}[theorem]{Lemma}
\newtheorem{proposition}[theorem]{Proposition}
\newtheorem{definition}[theorem]{Definition}
\newtheorem{corollary}[theorem]{Corollary}
\theoremstyle{remark}
\newtheorem{remark}[theorem]{Remark}
\title{Simple Yetter-Drinfeld Modules and Extensions
over Fractional Infinite-Dimensional Taft Algebras}
\author{Yiwei Zheng}

\date{}

\newcommand{\firstpageinfo}[1]{%
  \begingroup
  \renewcommand{\thefootnote}{}%
  \footnotetext{#1}%
  \addtocounter{footnote}{-1}%
  \endgroup
}

\begin{document}
\maketitle

\firstpageinfo{%
}

\begin{abstract}
We classify all simple Yetter-Drinfeld modules over a fractional
infinite-dimensional Taft Hopf algebra $T=T(\underline m,t,\xi)$,
including the infinite-dimensional induced family and the finite-dimensional
truncated and cyclic families. We determine their explicit structures and
isomorphism criteria, and classify first extensions between finite-dimensional
simple modules together with their middle modules. We also establish the primary decomposition with respect to $z^t$ and 
prove that every nonzero-eigenvalue component is wild when $m>1$.
\end{abstract}

\section{Introduction}

Let $\mathbbm{k}$ be an algebraically closed field of characteristic
zero. Yetter-Drinfeld modules play an important role in the theory
of Hopf algebras and braided tensor categories
\cite{Yetter1990,RadfordTowber1993}. They carry compatible
module and comodule structures, and their categories provide
a natural setting for braided Hopf algebras. In particular,
bosonization relates Hopf algebras in braided categories
to ordinary Hopf algebras \cite{Radford1985,Majid1994}.

Nichols algebras, originating in the work of Nichols
\cite{Nichols1978}, are fundamental objects in the structure
and classification theory of pointed Hopf algebras. They form
a central ingredient of the lifting method developed by
Andruskiewitsch and Schneider
\cite{AndruskiewitschSchneider2002,AndruskiewitschSchneider2010}.
Their structure has been studied extensively, including
the classification of arithmetic root systems associated
with Nichols algebras of diagonal type
\cite{Heckenberger2009}, Nichols algebras of semisimple
Yetter-Drinfeld modules
\cite{AndruskiewitschHeckenbergerSchneider2010},
and those arising from nonabelian groups
\cite{AndruskiewitschGrana2003}.

Taft algebras provide classical examples of finite-dimensional
pointed Hopf algebras \cite{Taft1971}. The representation
theory of their Drinfeld doubles and related finite-dimensional
Hopf algebras has been investigated in several works.
Chen studied finite-dimensional representations of a quantum
double associated with Taft algebras \cite{Chen2002},
while Erdmann, Green, Snashall and Taillefer studied
simple and indecomposable representations of Drinfeld doubles
of a related family of Hopf algebras
\cite{ErdmannEtAl2006}. Simple Yetter-Drinfeld modules
over finite-dimensional pointed Hopf algebras have also
been studied through cocycle twists
\cite{BenkartPereiraWitherspoon2010}.

In contrast, explicit classifications of simple
Yetter-Drinfeld modules over infinite-dimensional Hopf
algebras are less systematic. Recently, Zhen, Liu and Yu
classified simple Yetter-Drinfeld modules over
infinite-dimensional Taft algebras and investigated
their associated Nichols algebras
\cite{ZhenLiuYu2025}. They also obtained a classification
over generalized Liu algebras \cite{ZhenLiuYu2026}.
The fractional infinite-dimensional Taft Hopf algebras
considered in the present paper arise from the fraction
data appearing in Liu's study of prime Hopf algebras
of Gelfand--Kirillov dimension one \cite{Liu2020}.
Unlike the usual single-generator Taft setting, these
Hopf algebras possess several commuting skew-primitive
generators with a common power, leading to multivariable
coaction coefficients and additional possibilities
for simple modules and their extensions.

In this paper, we classify all simple left-left
Yetter-Drinfeld modules over the fractional
infinite-dimensional Taft Hopf algebra
$T=T(\underline m,t,\xi)$. The simple modules fall
into three families: infinite-dimensional induced
modules, finite-dimensional truncated modules,
and finite-dimensional cyclic modules. We determine
their explicit actions and coactions, dimensions,
standard elements, and isomorphism criteria.

We further compute the first extension spaces between
arbitrary finite-dimensional simple Yetter-Drinfeld
modules and describe the nonsplit middle modules
by generators, relations, bases, actions and coactions.
A generalized eigenspace decomposition for the
central operator $z^t$ separates extensions having
distinct central eigenvalues and reduces the problem
to the truncated and cyclic cases.

Finally, we investigate the representation-theoretic
complexity of the finite-dimensional Yetter-Drinfeld
category. Motivated by the theory of tame and wild
representation types \cite{Drozd1980,Ringel1984},
we construct, for $m>1$, a representation embedding
of $\operatorname{mod}_{\mathrm{fd}}\mathbbm{k}
\langle X,Y\rangle$ into every nonzero-eigenvalue
component. This proves wildness in the
representation-embedding sense.

All modules are left-left Yetter-Drinfeld modules unless
otherwise stated. Congruences in the parameters
$a,b,a',b'$ are modulo $n$, unless another modulus is
specified. We distinguish extension equivalence from
isomorphism of the middle modules.

\section{Preliminaries}

\subsection{Yetter-Drinfeld modules}
Let $H$ be a Hopf algebra with bijective antipode $S$. A left-left
Yetter-Drinfeld module over $H$ is simultaneously a left $H$-module and a
left $H$-comodule $V$, with $\delta(v)=v_{(-1)}\otimes v_{(0)}$, satisfying
\begin{equation}\label{eq:YD}
 \delta(h\cdot v)
 =h_{(1)}v_{(-1)}S(h_{(3)})\otimes h_{(2)}\cdot v_{(0)}
 \qquad(h\in H,\ v\in V).
\end{equation}
Write ${}_H^H\mathcal{YD}$ for the category of all such modules and
${}_H^H\mathcal{YD}_{\mathrm{fd}}$ for its finite-dimensional subcategory.

For a left $H$-module $V$, the tensor product
$H\otimes V$ carries the adjoint--diagonal action
\[
 h\triangleright(c\otimes v)
 =h_{(1)}cS(h_{(3)})\otimes h_{(2)}v.
\]
Thus the Yetter-Drinfeld condition \eqref{eq:YD}
is equivalent to the $H$-linearity of the coaction
$\delta:V\to H\otimes V$ with respect to this action.

\subsection{Fractional infinite-dimensional Taft Hopf algebra}
Let $m$ be a positive integer and let $m_1,\ldots,m_\theta$ be positive
integers. For $1\le i\le\theta$, set
\[
 e_i=\min\{a\in\mathbb N_{>0}:m\mid am_i\},\qquad
 A=\{(a_1,\ldots,a_\theta):0\le a_i<e_i\}.
\]

\begin{definition}\label{def:fraction}\cite[Definition~3.1]{Liu2020}
The tuple $\underline m=(m_1,\ldots,m_\theta)$ is called a
\emph{fraction of $m$ of length $\theta$} if the following hold:
\begin{enumerate}[label=\textup{(\arabic*)},leftmargin=2.5em]
\item $\gcd(e_i,m_i)=1$ for every $1\le i\le\theta$;
\item $\prod_{i=1}^{\theta}e_i=m$;
\item if $\mathbf a,\mathbf b\in A$ and $\mathbf a\ne\mathbf b$, then
\[
 \sum_{i=1}^{\theta}a_im_i\not\equiv
 \sum_{i=1}^{\theta}b_im_i\pmod m.
\]
\end{enumerate}
\end{definition}
Conditions
\textup{(2)}--\textup{(3)} imply the bijection
\begin{equation}\label{eq:unique}
 A\longrightarrow\mathbb Z/m\mathbb Z,\qquad
 (a_1,\ldots,a_\theta)\longmapsto\sum_{i=1}^{\theta}a_im_i\pmod m.
\end{equation}
In particular, the usual fraction properties give
\[
 m_0=\gcd(m_1,\ldots,m_\theta),\qquad
 \gcd(m,m_0)=1,\qquad e_i\mid m_j\quad(i\ne j),
\]
so $m\mid m_im_j/m_0$ for $i\ne j$.

Let $t$ be a positive integer and put $n=mt$. Choose a primitive $n$th root of unity
$\xi$ such that
\[
 \xi^{e_im_i/m_0}=\xi^{e_jm_j/m_0}
 \qquad(1\le i,j\le\theta).
\]
The Hopf algebra $T=T(\underline m,t,\xi)$ is generated by
$g,y_1,\ldots,y_\theta$, subject to
\begin{equation}\label{eq:T-relations}
 g^n=1,\quad y_ig=\xi^{m_i/m_0}gy_i,\quad
 y_iy_j=y_jy_i,\quad y_i^{e_i}=y_j^{e_j}.
\end{equation}
We write, throughout,
\[
 z=y_1^{e_1}=\cdots=y_\theta^{e_\theta}.
\]
The coalgebra and antipode are specified by
\[
\begin{gathered}
 \Delta(g)=g\otimes g,\qquad
 \Delta(y_i)=1\otimes y_i+y_i\otimes g^{tm_i},\\
 \varepsilon(g)=1,\quad\varepsilon(y_i)=0,\qquad
 S(g)=g^{-1},\quad S(y_i)=-y_i g^{-tm_i}.
\end{gathered}
\]
Since $\xi^{-tm_i^2/m_0}$ has order $e_i$, the quantum binomial
theorem applied to $\Delta(y_i)$ gives
$\Delta(z)=1\otimes z+z\otimes1$.
The PBW basis is
\begin{equation}\label{eq:PBWT}
 \{z^\ell y^{\mathbf r}g^j:\ell\ge0,\ 0\le r_i<e_i,
 \ 0\le j<n\},\qquad y^{\mathbf r}=\prod_i y_i^{r_i}.
\end{equation}
The positive grading with
$\deg(y_i)=m/e_i$ and $\deg(g)=0$ is compatible with the
coproduct; its degree-zero subcoalgebra is $\mathbbm{k}\langle g\rangle$.
It follows that $T$ is pointed and
$G(T)=\{g^a:a\in\mathbb Z/n\mathbb Z\}$.
We shall also use
$Z(T)=\mathbbm{k}[z^t]$ and
$\operatorname{ord}(\xi^{e_im_i/m_0})=t$;
see \cite[Proposition~4.1]{Liu2020}.

In particular,
\begin{equation}\label{eq:yYD}
\begin{aligned}
 \delta(y_iv)={}&-v_{(-1)}y_ig^{-tm_i}\otimes v_{(0)}
 +v_{(-1)}g^{-tm_i}\otimes y_iv_{(0)}\\
 &+y_iv_{(-1)}g^{-tm_i}\otimes g^{tm_i}v_{(0)}.
\end{aligned}
\end{equation}

\subsection{Gaussian binomial identities and PBW consequences}

For $q$ an indeterminate, let
\[
 [s]_q=1+q+\cdots+q^{s-1},\qquad
 [s]_q!=\prod_{j=1}^s[j]_q,
 \qquad
 {s\brack j}_q=\frac{[s]_q!}{[j]_q![s-j]_q!}.
\]
The Gaussian binomial coefficients belong to
$\mathbb Z[q]$ and satisfy
\begin{equation}\label{eq:Gaussian-recurrence}
\begin{aligned}
{s+1\brack j}_q
&={s\brack j}_q
 +q^{s+1-j}{s\brack j-1}_q,\\
{s+1\brack j}_q
&=q^j{s\brack j}_q
 +{s\brack j-1}_q.
\end{aligned}
\end{equation}
We use the conventions
${s\brack j}_q=0$ for $j<0$ or $j>s$ and
${s\brack0}_q={s\brack s}_q=1$.

In particular, since
$q_i=\xi^{-tm_i^2/m_0}$ has order $e_i$,
\[
 {e_i\brack j}_{q_i}=0
 \qquad(0<j<e_i).
\]
Since Gaussian binomial coefficients belong to
$\mathbb Z[q]$, the recurrence identities remain
valid after specialization to roots of unity,
even when the factorial quotients cannot be
evaluated directly.

\begin{lemma}\label{lem:PBW-free}
The algebra $T$ is a free module of rank $n^2$
over its center $\mathbbm{k}[z^t]$, with basis
\[
 \{z^\ell y^{\mathbf r}g^j:
 0\le\ell<t,\ 0\le r_i<e_i,\ 0\le j<n\}.
\]
Consequently, for every $\lambda\in\mathbbm{k}$
and $s\ge1$,
\[
 \dim_{\mathbbm{k}}
 T/(z^t-\lambda)^sT=sn^2.
\]
\end{lemma}

\begin{proof}
Writing each nonnegative integer uniquely as
$tq+\ell$, with $0\le\ell<t$, the PBW basis
\eqref{eq:PBWT} gives
\[
 T=\bigoplus_{\substack{0\le\ell<t\\
 0\le r_i<e_i\\0\le j<n}}
 \mathbbm{k}[z^t]z^\ell y^{\mathbf r}g^j.
\]
The rank is
$t(\prod_i e_i)n=tmn=n^2$.
Taking the quotient by $(z^t-\lambda)^s$
gives the dimension formula.
\end{proof}

\begin{lemma}\label{lem:PBW-weights}
The map
\[
\begin{aligned}
 \{0,\ldots,t-1\}\times
 \prod_{i=1}^{\theta}\{0,\ldots,e_i-1\}
 &\longrightarrow\mathbb Z/n\mathbb Z,\\
 (\ell,\mathbf r)
 &\longmapsto
 \frac{\ell e_1m_1+\sum_i r_im_i}{m_0}
 \pmod n
\end{aligned}
\]
is bijective.
\end{lemma}

\begin{proof}
Suppose two pairs $(\ell,\mathbf r)$ and
$(\ell',\mathbf r')$ have the same image.
Since $m\mid e_1m_1$, reduction modulo $m$
gives
\[
 \sum_i(r_i-r_i')m_i\equiv0\pmod m.
\]
By \eqref{eq:unique}, $\mathbf r=\mathbf r'$.
Hence
\[
 \xi^{(\ell-\ell')e_1m_1/m_0}=1.
\]
Since $\xi^{e_1m_1/m_0}$ has order $t$,
we obtain $\ell=\ell'$. Thus the map is
injective, and hence bijective, as its
domain has cardinality $t\prod_i e_i=n$.
\end{proof}

Since $z$ is primitive, \eqref{eq:YD} gives
\begin{equation}\label{eq:zYD}
 \delta(zv)
 =v_{(-1)}\otimes zv_{(0)}
 +(zv_{(-1)}-v_{(-1)}z)\otimes v_{(0)}.
\end{equation}
In particular, if a Yetter-Drinfeld module $W$
satisfies
\[
 \delta(W)\subseteq C_T(z)\otimes W,
\]
where $C_T(z)$ is the centralizer of $z$ in $T$,
then
\[
 \delta(z^tw)
 =(\operatorname{id}\otimes z^t)\delta(w).
\]
Since $z^t\in Z(T)$, it follows that
\[
 z^t\in
 \operatorname{End}_{{}_T^T\mathcal{YD}}(W).
\]

\subsection{Extensions and elementary facts}

For Yetter-Drinfeld modules $U,V$, we regard
$\operatorname{Ext}^1(U,V)$ as the space of
equivalence classes of short exact sequences
\[
 0\longrightarrow V\longrightarrow E
 \longrightarrow U\longrightarrow0,
\]
with addition given by the Baer sum.
Two extensions are equivalent if there is an
isomorphism between their middle modules
inducing the identity maps on $V$ and $U$.
The zero class is represented by a split extension.

\begin{lemma}\label{lem:extension-indecomp}
Let
\[
 0\longrightarrow S\longrightarrow E
 \longrightarrow S'\longrightarrow0
\]
be a nonsplit extension of simple Yetter-Drinfeld
modules. Then $E$ is indecomposable, and $S$
is its unique nonzero proper submodule.
\end{lemma}

\begin{proof}
Let $U$ be a nonzero proper submodule of $E$.
Since $E$ has composition length two, $U$ is
simple. If $U\ne S$, then $U\cap S=0$ and the
projection $U\to S'$ is an isomorphism, forcing
the extension to split. Hence $U=S$, and
$E$ is indecomposable.
\end{proof}

For finite-dimensional simple modules $S,S'$,
Schur's lemma also gives
\[
 \operatorname{Hom}(S,S')=
 \begin{cases}
 \mathbbm{k},&S\cong S',\\
 0,&S\not\cong S'.
 \end{cases}
\]

\section{Classification of simple Yetter-Drinfeld modules}

We classify all simple Yetter-Drinfeld modules over $T$ by means of
standard elements and quotients of induced modules.

\begin{definition}\cite[Definition~2.13]{ZhenLiuYu2025}
A nonzero element $v_0$ of a Yetter-Drinfeld module is called a
\emph{standard element of type $(a,b)$} if
\[
 \delta(v_0)=g^a\otimes v_0,\qquad gv_0=\xi^b v_0,
 \qquad a,b\in\mathbb Z/n\mathbb Z.
\]

\end{definition}

\begin{lemma}\label{lem:standard-exists}
Every nonzero Yetter-Drinfeld module $V$ over $T$ contains a standard
element. If $V$ is simple, then $V=Tv_0$ for some standard element
$0\ne v_0\in V$.
\end{lemma}
\begin{proof}
Since $T$ is pointed, $V$ contains a nonzero
one-dimensional subcomodule. Thus for some $a\in\mathbb Z/n\mathbb Z$ and $0\ne v\in V$,
\[
 \delta(v)=g^a\otimes v.
\]
The Yetter-Drinfeld condition yields
\[
 \delta(g^jv)=g^a\otimes g^jv\qquad(j\ge0).
\]
Since $g^n=1$ and $X^n-1$ has distinct roots in $\mathbbm{k}$, we can choose
$0\ne v_0\in\operatorname{span}_{\mathbbm{k}}\{v,gv,\ldots,g^{n-1}v\}$ such that
\[
 \delta(v_0)=g^a\otimes v_0,\qquad gv_0=\xi^b v_0
\]
for some $b\in\mathbb Z/n\mathbb Z$. Hence $v_0$ is standard. Moreover,
\[
 \delta(xv_0)=x_{(1)}g^aS(x_{(3)})\otimes x_{(2)}v_0
 \in T\otimes Tv_0\qquad(x\in T).
\]
Thus $Tv_0$ is a nonzero Yetter-Drinfeld submodule, and $V=Tv_0$
when $V$ is simple.
\end{proof}

For $a,b\in\mathbb Z/n\mathbb Z$, let $\mathbbm{k}v_0$ be the one-dimensional
$\mathbbm{k}\langle g\rangle$-module and comodule defined by
$gv_0=\xi^bv_0$ and $\delta(v_0)=g^a\otimes v_0$. Put
\begin{equation}\label{eq:M}
 M(a,b)=T\otimes_{\mathbbm{k}\langle g\rangle}\mathbbm{k}v_0,\qquad v=1\otimes v_0.
\end{equation}
For $x,u\in T$, define the action and coaction by
\[
 x\cdot(uv)=(xu)v,\qquad
 \delta(xv)=x_{(1)}g^aS(x_{(3)})\otimes x_{(2)}v.
\]
These formulas are well defined. Indeed, the defining relation
$ugv=\xi^buv$ is respected by the coaction, since
\[
\begin{aligned}
 \delta(ugv)
 &=u_{(1)}gg^ag^{-1}S(u_{(3)})\otimes u_{(2)}gv\\
 &=\xi^b u_{(1)}g^aS(u_{(3)})\otimes u_{(2)}v
 =\xi^b\delta(uv).
\end{aligned}
\]
The Hopf identities give
\[
\begin{aligned}
 (\Delta\otimes\operatorname{id})\delta(xv)
 &=(\operatorname{id}\otimes\delta)\delta(xv)\\
 &=x_{(1)}g^aS(x_{(5)})\otimes
   x_{(2)}g^aS(x_{(4)})\otimes x_{(3)}v,\\
 (\varepsilon\otimes\operatorname{id})\delta(xv)&=xv.
\end{aligned}
\]
The displayed coaction is $T$-linear for the adjoint diagonal action
from Section~2.1. Hence $M(a,b)$ is a Yetter-Drinfeld module.
By the PBW basis of $T$,
\begin{equation}\label{eq:PBWM}
 \{z^\ell y^{\mathbf r}v:\ell\ge0,\ 0\le r_i<e_i\}
\end{equation}
is a basis of $M(a,b)$. For any simple module containing a standard
$v_0$ of type $(a,b)$, the assignment $xv\mapsto xv_0$ defines a
surjective Yetter-Drinfeld morphism $M(a,b)\twoheadrightarrow V$.

\begin{lemma}\label{lem:C}
For $\mathbf s\in\mathbb N^\theta$,
\begin{equation}\label{eq:C-coaction}
 \delta(y^{\mathbf s}v)=\sum_{\mathbf0\le\mathbf j\le\mathbf s}
 C_{\mathbf s,\mathbf j}(a,b)
 y^{\mathbf j}g^{a-t\sum_i s_im_i}\otimes y^{\mathbf s-\mathbf j}v,
\end{equation}
where
\begin{equation}\label{eq:C}
 C_{\mathbf s,\mathbf j}(a,b)=\prod_{i=1}^\theta\left[
 {s_i\brack j_i}_{\xi^{-tm_i^2/m_0}}
 \prod_{k=0}^{j_i-1}
 \left(\xi^{btm_i}-
 \xi^{-am_i/m_0+tm_i^2(s_i-j_i+k)/m_0}\right)\right].
\end{equation}
For a fixed $i$, let $C^{(i)}_{s,j}(a,b)$ denote the
one-variable factor in brackets in \eqref{eq:C}. The corresponding
one-variable statement is
\[
 \delta(y_i^s v)=\sum_{j=0}^{s}C^{(i)}_{s,j}(a,b)
 y_i^j g^{a-stm_i}\otimes y_i^{s-j}v,
 \qquad C^{(i)}_{0,0}=1.
\]
Here $\mathbf1_i$ is the $i$th unit multi-index, and Gaussian
binomial coefficients are evaluated as polynomials.
\end{lemma}

\begin{proof}
We first prove the one-variable formula.
Fix $i$. For $s=0$, we have
\[
\delta(v)=g^a\otimes v,
\qquad C^{(i)}_{0,0}(a,b)=1.
\]
For $s=1$, \eqref{eq:yYD} gives
\[
\begin{aligned}
\delta(y_iv)
={}&g^{a-tm_i}\otimes y_iv+\left(\xi^{btm_i}-\xi^{-am_i/m_0}\right)
y_i g^{a-tm_i}\otimes v.
\end{aligned}
\]

Suppose that
\[
\delta(y_i^sv)
=\sum_{j=0}^{s}C^{(i)}_{s,j}(a,b)
y_i^jg^{a-stm_i}\otimes y_i^{s-j}v.
\]
Applying \eqref{eq:yYD}, we obtain
\[
\begin{aligned}
\delta(y_i^{s+1}v)
={}&\sum_{j=0}^{s}C^{(i)}_{s,j}(a,b)
y_i^jg^{a-(s+1)tm_i}
\otimes y_i^{s+1-j}v\\
&+\sum_{j=0}^{s}C^{(i)}_{s,j}(a,b)
\left(
\xi^{btm_i-tm_i^2(s-j)/m_0}
-\xi^{-am_i/m_0+tm_i^2s/m_0}
\right)
y_i^{j+1}g^{a-(s+1)tm_i}
\otimes y_i^{s-j}v.
\end{aligned}
\]
Consequently,
\[
\begin{aligned}
C^{(i)}_{s+1,j}(a,b)
={}&C^{(i)}_{s,j}(a,b)+\left(
\xi^{btm_i-tm_i^2(s-j+1)/m_0}
-\xi^{-am_i/m_0+tm_i^2s/m_0}
\right)
C^{(i)}_{s,j-1}(a,b),
\end{aligned}
\]
where $C^{(i)}_{s,j}=0$ if $j<0$ or
$j>s$.

Substituting \eqref{eq:C} into the recurrence
and applying the Gaussian binomial identities
in \eqref{eq:Gaussian-recurrence}, we obtain
the required product formula. The boundary
cases follow directly. Since these identities
are polynomial identities, the formula
remains valid for $s\ge e_i$.

We now consider the multivariable case.
For $i\ne j$, the fraction conditions imply
\[
\xi^{tm_im_j/m_0}=1,
\qquad
g^{tm_i}y_j=y_jg^{tm_i}.
\]
Assume that \eqref{eq:C-coaction} holds
for $\mathbf s$. Applying \eqref{eq:yYD}
to $y_i y^{\mathbf s}v$ gives
\[
\begin{aligned}
\delta(y_i y^{\mathbf s}v)
={}&\sum_{\mathbf0\le\mathbf j\le\mathbf s}
C_{\mathbf s,\mathbf j}(a,b)
y^{\mathbf j}
g^{a-t(\sum_k s_km_k+m_i)}
\otimes
y^{\mathbf s+\mathbf1_i-\mathbf j}v\\
&+\sum_{\mathbf0\le\mathbf j\le\mathbf s}
C_{\mathbf s,\mathbf j}(a,b)
\left(
\xi^{btm_i-tm_i^2(s_i-j_i)/m_0}
-\xi^{-am_i/m_0+tm_i^2s_i/m_0}
\right)\\
&\qquad\cdot
y^{\mathbf j+\mathbf1_i}
g^{a-t(\sum_k s_km_k+m_i)}
\otimes y^{\mathbf s-\mathbf j}v.
\end{aligned}
\]
Comparing coefficients yields precisely
the one-variable recurrence in the $i$th
coordinate, while the remaining factors
are unchanged. Thus induction on
$\sum_i s_i$ proves \eqref{eq:C-coaction}.
\end{proof}

\begin{lemma}\label{lem:di}
If $t\nmid a$, then $C_{\mathbf s,\mathbf s}(a,b)\ne0$ for every
$\mathbf s$. If $t\mid a$, for each $i$ there is a unique
$1\le d_i\le e_i$ satisfying
\begin{equation}\label{eq:di}
 \frac at+bm_0\equiv m_i(d_i-1)\pmod{e_i}.
\end{equation}
Moreover,
\begin{equation}\label{eq:di-properties}
 C^{(i)}_{s,s}(a,b)\ne0\ (1\le s<d_i),\qquad
 C^{(i)}_{d_i,j}(a,b)=0\ (1\le j\le d_i),
\end{equation}
and $y_i^{d_i}v$ is a nonzero standard element of type
$(a-td_im_i,b-d_im_i/m_0)$ in $M(a,b)$.
\end{lemma}
\begin{proof}
By \eqref{eq:C},
\[
 C^{(i)}_{s,s}(a,b)=\prod_{j=0}^{s-1}
 \left(\xi^{btm_i}-\xi^{-am_i/m_0+tm_i^2j/m_0}\right).
\]
The first possible zero occurs exactly when
\[
 \xi^{btm_i}=\xi^{-am_i/m_0+tm_i^2(d_i-1)/m_0}.
\]
Since $\xi^{tm_i^2/m_0}$ has order $e_i$, this equality has a
solution modulo $e_i$ if and only if
$\xi^{ae_im_i/m_0}=1$, equivalently $t\mid a$. If $t\mid a$,
the congruence reduces to \eqref{eq:di}, which has exactly one
solution $d_i\in\{1,\ldots,e_i\}$ because
$\gcd(m_i,e_i)=1$. In \eqref{eq:C} with $s_i=d_i$ and
$1\le j_i\le d_i$, the factor indexed by $k=j_i-1$ is zero.
Consequently,
\[
 \delta(y_i^{d_i}v)=g^{a-td_im_i}\otimes y_i^{d_i}v,
 \qquad gy_i^{d_i}v=\xi^{b-d_im_i/m_0}y_i^{d_i}v.
\]
Nonvanishing follows from \eqref{eq:PBWM}, including the case
$d_i=e_i$, when $y_i^{d_i}v=zv$.
\end{proof}

\begin{lemma}\label{lem:z-simple}
Let $V=Tv_0$ be simple with $v_0$ standard of type $(a,b)$ and
$t\mid a$. Then $zV$ is a Yetter-Drinfeld submodule, $z^t$ acts
as some scalar $\lambda\in\mathbbm{k}$, and $\dim_{\mathbbm{k}} V\le n$. The identity
$\delta(zw)=(\operatorname{id}\otimes z)\delta(w)$ also holds
on $M(a,b)$ when $t\mid a$.
\end{lemma}

\begin{proof}
Since $t\mid a$, the element $z$ commutes
with every first coaction factor in
\eqref{eq:C-coaction}. As $z$ is primitive,
the Yetter-Drinfeld compatibility gives
\[
\delta(zw)=(\operatorname{id}\otimes z)\delta(w)
\]
for $w\in V$ or $M(a,b)$.
Since $zT=Tz$, it follows that $zV$
is a Yetter-Drinfeld submodule.

Moreover, $z^t\in Z(T)$, and the above
identity shows that its action is a
Yetter-Drinfeld endomorphism.

By \eqref{eq:PBWM}, $M(a,b)$ is a free
$\mathbbm{k}[z^t]$-module of rank
$t\prod_i e_i=n$. Hence $V$ is a
nonzero finitely generated
$\mathbbm{k}[z^t]$-module.
Choose a maximal proper
$\mathbbm{k}[z^t]$-submodule $N$ of $V$.
Since $\mathbbm{k}$ is algebraically
closed, there exists
$\lambda\in\mathbbm{k}$ such that
\[
V/N\cong
\mathbbm{k}[z^t]/(z^t-\lambda).
\]
Consequently,
\[
(z^t-\lambda)V\subseteq N\ne V.
\]
Since $z^t-\lambda$ acts as a
Yetter-Drinfeld endomorphism, its image
is a proper Yetter-Drinfeld submodule.
The simplicity of $V$ therefore implies
\[
(z^t-\lambda)V=0.
\]
Thus $V$ is a quotient of
$M(a,b)/(z^t-\lambda)M(a,b)$, and
\[
\dim_{\mathbbm{k}}V
\le t\prod_i e_i=n.
\]
\end{proof}

\begin{proposition}\label{prop:simple-quotients}
For $t\mid a$, with $d_i$ from \eqref{eq:di}, define
\begin{equation}\label{eq:V0Vlambda}
 V_0(a,b)=M(a,b)\Big/\sum_iT(y_i^{d_i}v),\qquad
 V_\lambda(a,b)=M(a,b)/(z^t-\lambda)M(a,b)
 \quad(\lambda\in\mathbbm{k}^\times).
\end{equation}
Both are Yetter-Drinfeld modules, with respective bases
\begin{equation}\label{eq:finite-bases}
 \{y^{\mathbf r}\bar v:0\le r_i<d_i\},\qquad
 \{z^\ell y^{\mathbf r}\bar v:0\le\ell<t,\ 0\le r_i<e_i\},
\end{equation}
and dimensions $\prod_i d_i$ and $n$, respectively.
\end{proposition}

\begin{proof}
By Lemma~\ref{lem:di}, each
$y_i^{d_i}v$ is standard, so
$T(y_i^{d_i}v)$ is a Yetter-Drinfeld
submodule. Since $y_i^{e_i}=z$
and $d_i\le e_i$, we have
$zv\in T(y_i^{d_i}v)$. The PBW basis
therefore gives
\[
\sum_iT(y_i^{d_i}v)
=\operatorname{span}_{\mathbbm{k}}
\left\{
z^\ell y^{\mathbf r}v:
\ell\ge1\ \text{or}\
r_i\ge d_i\ \text{for some }i
\right\}.
\]
Thus $V_0(a,b)$ has the stated basis.

For $\lambda\in\mathbbm{k}^{\times}$,
Lemma~\ref{lem:z-simple} shows that
$(z^t-\lambda)M(a,b)$ is a
Yetter-Drinfeld submodule.
Since
\[
M(a,b)=
\bigoplus_{\substack{0\le\ell<t\\
0\le r_i<e_i}}
\mathbbm{k}[z^t]z^\ell y^{\mathbf r}v,
\]
the quotient $V_\lambda(a,b)$ has
the stated basis. In particular,
\[
\dim_{\mathbbm{k}}V_0(a,b)
=\prod_i d_i,\qquad
\dim_{\mathbbm{k}}V_\lambda(a,b)=n.
\]
\end{proof}

\begin{theorem}\label{thm:simples}
Every simple Yetter-Drinfeld $T$-module belongs to one of
the following three families:
\begin{equation}\label{eq:three-families}
 \begin{array}{ll}
 M(a,b),&t\nmid a,\qquad \dim_{\mathbbm{k}}M(a,b)=\infty,\\
 V_0(a,b),&t\mid a,\qquad \dim_{\mathbbm{k}}V_0(a,b)=\prod_i d_i,\\
 V_\lambda(a,b),&t\mid a,\ \lambda\in\mathbbm{k}^\times,\quad
 \dim_{\mathbbm{k}}V_\lambda(a,b)=n.
 \end{array}
\end{equation}
The action of $g$ on the finite families is
\begin{equation}\label{eq:g-on-basis}
 g(z^\ell y^{\mathbf r}\bar v)
 =\xi^{b-(\ell e_1m_1+\sum_i r_im_i)/m_0}
 z^\ell y^{\mathbf r}\bar v,
\end{equation}
where $\ell=0$ for $V_0$. The actions of $y_i$ are
\begin{equation}\label{eq:simple-action}
\begin{aligned}
 y_i(y^{\mathbf r}\bar v)&=\begin{cases}
 y^{\mathbf r+\mathbf1_i}\bar v,&r_i<d_i-1,\\
 0,&r_i=d_i-1
 \end{cases}\quad\text{on }V_0,\\
 y_i(z^\ell y^{\mathbf r}\bar v)&=\begin{cases}
 z^\ell y^{\mathbf r+\mathbf1_i}\bar v,&r_i<e_i-1,\\
 z^{\ell+1}\prod_{j\ne i}y_j^{r_j}\bar v,
 &r_i=e_i-1,\ \ell<t-1,\\
 \lambda\prod_{j\ne i}y_j^{r_j}\bar v,
 &r_i=e_i-1,\ \ell=t-1
 \end{cases}\quad\text{on }V_\lambda.
\end{aligned}
\end{equation}
Their coactions follow from \eqref{eq:C-coaction} and
$\delta(z^\ell y^{\mathbf r}\bar v)
 =(\operatorname{id}\otimes z^\ell)\delta(y^{\mathbf r}\bar v)$.
For $M(a,b)$ the action is induced from $T$ and the coaction is
\eqref{eq:C-coaction}.
\end{theorem}

\begin{proof}
We first prove simplicity.

Suppose $t\nmid a$, and let
$0\ne U\subseteq M(a,b)$ be a
Yetter-Drinfeld submodule.
By Lemma~\ref{lem:standard-exists},
there exists $0\ne w\in U$ such that
$\delta(w)=g^j\otimes w$.
Write
\[
w=\sum_{\ell,\mathbf r}
u_{\ell,\mathbf r}z^\ell y^{\mathbf r}v,
\qquad
\mathbf s=(e_1\ell+r_1,r_2,\ldots,r_\theta).
\]
Comparing the coefficient of $v$
in the second tensor factor gives
\[
\sum_{\ell,\mathbf r}
u_{\ell,\mathbf r}C_{\mathbf s,\mathbf s}(a,b)
z^\ell y^{\mathbf r}
g^{a-t\sum_i r_im_i}
=u_{0,\mathbf0}g^j,
\]
where we used $n\mid te_1m_1$.
By Lemma~\ref{lem:di},
$C_{\mathbf s,\mathbf s}(a,b)\ne0$.
The PBW monomials on the left are
pairwise distinct, so
$u_{\ell,\mathbf r}=0$ unless
$(\ell,\mathbf r)=(0,\mathbf0)$.
Hence $w\in\mathbbm{k}^{\times}v$,
and $U=M(a,b)$.

We next consider $V_0(a,b)$ and $V_\lambda(a,b)$.
Let $U$ be a nonzero Yetter--Drinfeld submodule of either module.
By Lemma~\ref{lem:PBW-weights}, the reduced PBW basis vectors
have pairwise distinct $g$-weights, so $U$ contains one of them.

For $V_0(a,b)$, let
$y^{\mathbf r}\bar v\in U$.
Since $\mathbf r<\mathbf d$,
Lemma~\ref{lem:di} gives
$C_{\mathbf r,\mathbf r}(a,b)\ne0$.
The coaction of this vector contains
\[
C_{\mathbf r,\mathbf r}(a,b)
y^{\mathbf r}g^{a-t\sum_i r_im_i}
\otimes\bar v.
\]
As $U$ is a subcomodule, comparison
of PBW coefficients gives
$\bar v\in U$. Hence $U=V_0(a,b)$.

For $V_\lambda(a,b)$, we have
\[
y_i^{e_it}=z^t=\lambda\ne0.
\]
Thus every $y_i$ acts invertibly,
with inverse given by a polynomial
in $y_i$. Every reduced PBW vector
therefore generates $\bar v$
under the $T$-action. Hence
$U=V_\lambda(a,b)$.

It remains to prove exhaustiveness.
Let $V$ be simple and choose a
standard element $v_0$ of type
$(a,b)$. The induced morphism
$M(a,b)\twoheadrightarrow V$
is surjective.

If $t\nmid a$, the simplicity of
$M(a,b)$ makes this morphism
an isomorphism.

Suppose $t\mid a$. By
Lemma~\ref{lem:z-simple}, $z^t$
acts on $V$ as a scalar
$\lambda\in\mathbbm{k}$.
If $\lambda\ne0$, the morphism
factors through $V_\lambda(a,b)$,
which is simple, and hence
$V\cong V_\lambda(a,b)$.

If $\lambda=0$, the submodule
$zV$ is either zero or $V$.
The latter would imply
$z^tV=V$, a contradiction.
Thus $zV=0$.

By Lemma~\ref{lem:di},
$y_i^{d_i}v_0$ is a standard
element whenever it is nonzero.
If $y_i^{d_i}v_0\ne0$, simplicity
and normality of $y_i^{d_i}$
would give
\[
V=T(y_i^{d_i}v_0)
=y_i^{d_i}V.
\]
But $y_i^{e_i}V=zV=0$, so
$y_i^{d_i}$ is nilpotent on $V$
and cannot be surjective.
Therefore
\[
y_i^{d_i}v_0=0
\qquad(1\le i\le\theta).
\]
The morphism from $M(a,b)$
factors through $V_0(a,b)$,
which is simple. Hence
$V\cong V_0(a,b)$.
\end{proof}

\begin{proposition}\label{prop:iso}
The modules $M(a,b)$ and $V_0(a,b)$ are isomorphic to another
member of the same family only when their two parameters agree
modulo $n$. Moreover,
\begin{equation}\label{eq:cyclic-iso}
\begin{aligned}
 V_\lambda(a,b)\cong V_\mu(a',b')\iff{}&\lambda=\mu
 \text{ and for some }0\le\ell<t,\ J\subseteq\{i:d_i<e_i\},\\
 a'&\equiv a-t\sum_{i\in J}d_im_i\pmod n,\\
 b'&\equiv b-\frac{\ell e_1m_1+\sum_{i\in J}d_im_i}{m_0}\pmod n.
\end{aligned}
\end{equation}
\end{proposition}

\begin{proof}
For the infinite family $M(a,b)$
with $t\nmid a$, the proof of
Theorem~\ref{thm:simples} shows
that every standard element is
a scalar multiple of $v$.

The same holds for $V_0(a,b)$:
its reduced PBW vectors have
distinct $g$-weights, and
$C_{\mathbf r,\mathbf r}(a,b)\ne0$
for $\mathbf0\le\mathbf r<\mathbf d$.
Thus no nonconstant reduced PBW
vector is coaction-homogeneous.

Consequently, an isomorphism
between two members of either
family must preserve the type
$(a,b)$. The parameters therefore
agree modulo $n$.

Now suppose
\[
f:V_\mu(a',b')
\xrightarrow{\sim}V_\lambda(a,b).
\]
Since $f(\bar v')$ is standard
of type $(a',b')$ and the reduced
PBW vectors have distinct
$g$-weights, we have
\[
f(\bar v')
=cz^\ell y^{\mathbf r}\bar v,
\qquad c\in\mathbbm{k}^{\times},
\quad 0\le\ell<t,\quad
0\le r_i<e_i.
\]

By \eqref{eq:C-coaction},
this vector is coaction-homogeneous
if and only if
\[
C^{(i)}_{r_i,1}(a,b)=0
\qquad\text{for every }r_i>0.
\]
Indeed, for $0<r_i<e_i$, the
Gaussian integer
$[r_i]_{\xi^{-tm_i^2/m_0}}$
is nonzero. Lemma~\ref{lem:di}
therefore gives
\[
r_i=0
\quad\text{or}\quad
r_i=d_i<e_i.
\]
Conversely, these conditions
make all positive coaction
coefficients vanish by
\eqref{eq:di-properties}.

Hence
\[
f(\bar v')
=cz^\ell\prod_{i\in J}
y_i^{d_i}\bar v
\]
for some
$J\subseteq\{i:d_i<e_i\}$.
Comparing its coaction degree
and $g$-weight gives
\[
\begin{aligned}
a'&\equiv
a-t\sum_{i\in J}d_im_i\pmod n,\\
b'&\equiv
b-\frac{\ell e_1m_1+
\sum_{i\in J}d_im_i}{m_0}\pmod n.
\end{aligned}
\]
Moreover, since $z^t$ acts as
$\lambda$ and $\mu$ on the two
modules, respectively, an
isomorphism requires $\lambda=\mu$.
Thus \eqref{eq:cyclic-iso} is necessary.

Conversely, suppose the conditions
in \eqref{eq:cyclic-iso} hold.
Then
\[
w=z^\ell\prod_{i\in J}
y_i^{d_i}\bar v
\]
is a nonzero standard element
of type $(a',b')$.

By the universal property of
$M(a',b')$, the assignment
$v'\mapsto w$ induces a nonzero
Yetter-Drinfeld morphism
\[
M(a',b')\longrightarrow V_\lambda(a,b).
\]
Since $\lambda=\mu$ and
$(z^t-\mu)w=0$, it factors through
$V_\mu(a',b')$.
We obtain a nonzero morphism
\[
V_\mu(a',b')
\longrightarrow V_\lambda(a,b).
\]
Both modules are simple, so this
morphism is an isomorphism.
\end{proof}

\section{First extensions of finite-dimensional simple modules}

We determine first extensions between finite-dimensional simple modules.
The central endomorphism $z^t$ first separates the eigenvalue
components. Then we describe the extensions within the truncated and
cyclic families, including the isomorphism types of their middle modules.

\subsection{Generalized eigenspace decomposition}

For $0\le r<t$, set
\begin{equation}\label{eq:Tr}
 T^{(r)}=\operatorname{span}_{\mathbbm{k}}
 \{z^\ell y^{\mathbf s}g^j:\ell\ge0,\ 0\le s_i<e_i,\
 j\equiv r\pmod t\}.
\end{equation}

\begin{lemma}\label{lem:coalgebra-splitting}
The vector-space decomposition $T=\bigoplus_{r=0}^{t-1}T^{(r)}$
is a direct sum of subcoalgebras. Every left $T$-comodule $W$ has a
functorial decomposition $W=\bigoplus_r W^{(r)}$ with
$\delta(W^{(r)})\subseteq T^{(r)}\otimes W^{(r)}$. If $W$ is a
Yetter-Drinfeld module, all $W^{(r)}$ are Yetter-Drinfeld submodules.
\end{lemma}
\begin{proof}
For $z^\ell y^{\mathbf s}g^j\in T^{(r)}$, the coproducts
\[
 \Delta(y_ig^j)=g^j\otimes y_ig^j+y_ig^j\otimes g^{j+tm_i},
 \qquad
 \Delta(zg^j)=g^j\otimes zg^j+zg^j\otimes g^j
\]
show, using multiplicativity of $\Delta$, that
\[
 \Delta(T^{(r)})\subseteq T^{(r)}\otimes T^{(r)},
 \qquad T=\bigoplus_{r=0}^{t-1}T^{(r)}.
\]
For a left comodule $W$, let $p_r:T\to \mathbbm{k}$ be $\varepsilon$
on $T^{(r)}$ and zero on the other summands, and put
$P_r=(p_r\otimes\operatorname{id})\delta$. Then
\[
 P_rP_s=\delta_{r,s}P_r,\qquad
 \sum_{r=0}^{t-1}P_r=\operatorname{id}_W,\qquad
 W^{(r)}=P_rW.
\]
Coassociativity gives
$\delta(W^{(r)})\subseteq T^{(r)}\otimes W^{(r)}$.
Since $gT^{(r)}g^{-1}=T^{(r)}$ and each $g^{\pm tm_i}$
has group exponent divisible by $t$, the formulas
\eqref{eq:YD} and \eqref{eq:yYD} yield
\[
 gW^{(r)}\subseteq W^{(r)},\qquad
 y_iW^{(r)}\subseteq W^{(r)}\quad(1\le i\le\theta).
\]
Thus each $W^{(r)}$ is a Yetter-Drinfeld submodule.
\end{proof}

\begin{theorem}\label{thm:primary}
For every finite-dimensional Yetter-Drinfeld module $W$,
\begin{equation}\label{eq:W0}
 W=W^{(0)},\qquad z^t\in\operatorname{End}_{{}_T^T\mathcal{YD}}(W).
\end{equation}
If $\lambda_1,\ldots,\lambda_s$ are the eigenvalues of $z^t$ and
$N\ge\dim_{\mathbbm{k}}W$, then
\begin{equation}\label{eq:primary}
 W=\bigoplus_{j=1}^s\ker(z^t-\lambda_j)^N,
\end{equation}
where every summand is a Yetter-Drinfeld submodule.
\end{theorem}
\begin{proof}
By Lemma~\ref{lem:coalgebra-splitting},
$W=\bigoplus_{r=0}^{t-1}W^{(r)}$ as Yetter-Drinfeld modules.
If $W^{(r)}\ne0$ for $r\ne0$, finite dimensionality gives a
simple submodule $S\subseteq W^{(r)}$ with a standard vector
$v_0$ of some type $(a,b)$. Then
\[
 \delta(v_0)=g^a\otimes v_0\in T^{(r)}\otimes W^{(r)}
 \ \Longrightarrow\ a\equiv r\not\equiv0\pmod t.
\]
Theorem~\ref{thm:simples} would imply $\dim_{\mathbbm{k}}S=\infty$, a
contradiction. Therefore $W=W^{(0)}$.

For $x\in T^{(0)}$,
\[
 zy_i=y_iz,\quad zg^t=\xi^{te_1m_1/m_0}g^tz=g^tz
 \ \Longrightarrow\ zx=xz.
\]
Since $\delta(W)\subseteq T^{(0)}\otimes W$ and $z$ is
primitive, the Yetter-Drinfeld identity gives
\[
 \delta(zw)=w_{(-1)}\otimes zw_{(0)}+
 (zw_{(-1)}-w_{(-1)}z)\otimes w_{(0)}
 =(\operatorname{id}\otimes z)\delta(w).
\]
Thus $z^t$ is a comodule endomorphism and, because
$z^t\in Z(T)$, also a module endomorphism.
The primary decomposition of the finite-dimensional operator
$z^t$ yields \eqref{eq:primary}. Each
$(z^t-\lambda_j)^N$ is a Yetter-Drinfeld endomorphism, so its
kernel is a Yetter-Drinfeld submodule.
\end{proof}

\begin{corollary}\label{cor:unequal}
If finite-dimensional simple Yetter-Drinfeld modules $S,S'$ have
distinct $z^t$-eigenvalues, then $\operatorname{Ext}^1(S',S)=0$.
In particular, all truncated--cyclic extensions split in both
directions.
\end{corollary}
\begin{proof}
In $0\to S\to E\to S'\to0$ with eigenvalues $\lambda\ne\mu$,
$(z^t-\lambda)(z^t-\mu)E=0$. The two factors are coprime, so
$E=\ker(z^t-\lambda)\oplus\ker(z^t-\mu)$, a direct sum of
Yetter-Drinfeld submodules mapping to the two endpoints.
\end{proof}

\subsection{Extensions between truncated simple modules}

Consider
\begin{equation}\label{eq:A-extension}
 0\longrightarrow V_0(a,b)=Tw_0\longrightarrow E
 \xrightarrow{p}V_0(a',b')=Tp(w_1)\longrightarrow0.
\end{equation}
Write $d_i'=d_i(a',b')$, choose $gw_1=\xi^{b'}w_1$,
and use $\mathbf r<\mathbf d$ to mean $0\le r_i<d_i$ for all $i$.
Its basis, obtained by lifting the quotient basis, is
\begin{equation}\label{eq:A-basis}
 \{y^{\mathbf r}w_0:0\le r_i<d_i\}\ \cup\
 \{y^{\mathbf s}w_1:0\le s_i<d_i'\}.
\end{equation}
Initially the additional action and coaction data are
\begin{equation}\label{eq:A-data}
 y_i^{d_i'}w_1=\sum_{\mathbf r<\mathbf d}u_{i,\mathbf r}
 y^{\mathbf r}w_0,\qquad
 \delta(w_1)=g^{a'}\otimes w_1+
 \sum_{\mathbf r<\mathbf d}c_{\mathbf r}\otimes y^{\mathbf r}w_0.
\end{equation}
The counit, coassociativity, and $g$-compatibility give
\begin{equation}\label{eq:A-c-eq}
\begin{aligned}
 \varepsilon(c_{\mathbf r})&=0,\\
 \Delta(c_{\mathbf r})&=g^{a'}\otimes c_{\mathbf r}
 +\sum_{\mathbf r\le\mathbf s<\mathbf d}
 C_{\mathbf s,\mathbf s-\mathbf r}(a,b)c_{\mathbf s}
 \otimes y^{\mathbf s-\mathbf r}g^{a-t\sum_i s_im_i},\\
 gc_{\mathbf r}g^{-1}&=
 \xi^{b'-b+\sum_i r_im_i/m_0}c_{\mathbf r}.
\end{aligned}
\end{equation}
For $c\in T$ and $v\in V_0(a,b)$ write
\begin{equation}\label{eq:D}
 \mathcal D_i(c\otimes v)=-cy_ig^{-tm_i}\otimes v
 +cg^{-tm_i}\otimes y_iv+y_icg^{-tm_i}\otimes g^{tm_i}v.
\end{equation}
Induction in \eqref{eq:yYD} gives the boundary equation
\begin{equation}\label{eq:A-master}
 \sum_{\mathbf r<\mathbf d}\mathcal D_i^{d_i'}
 (c_{\mathbf r}\otimes y^{\mathbf r}w_0)
 =\sum_{\mathbf r<\mathbf d}u_{i,\mathbf r}
 \left(\delta(y^{\mathbf r}w_0)
 -g^{a'-td_i'm_i}\otimes y^{\mathbf r}w_0\right).
\end{equation}

\begin{lemma}\label{lem:A-normal}
Assume $m>1$. A lift $w_1$ of the quotient standard generator
can be chosen so that
\begin{equation}\label{eq:A-normal}
 \delta(w_1)=g^{a'}\otimes w_1+\gamma zg^a\otimes w_0,
 \qquad y_i^{d_i'}w_1=u_iw_0\quad(1\le i\le\theta).
\end{equation}
Here $\gamma,u_i\in\mathbbm{k}$. If $\gamma\ne0$, then
\begin{equation}\label{eq:A-gamma}
 a'=a,\qquad b'-b\equiv-\frac{e_1m_1}{m_0}\pmod n,
 \qquad d_i'=d_i\quad(1\le i\le\theta).
\end{equation}
If $u_i\ne0$, then
\begin{equation}\label{eq:A-u}
 a\equiv a'-td_i'm_i,\qquad
 b\equiv b'-\frac{d_i'm_i}{m_0}\pmod n.
\end{equation}
\end{lemma}

\begin{proof}
Let $\pi:T\to\mathbbm{k}\langle g\rangle$
be the Hopf projection with $\pi(y_i)=0$.
By the projected coaction decomposition
and the diagonalizability of $g$, we may
choose $w_1$ such that
\[
gw_1=\xi^{b'}w_1,\qquad
(\pi\otimes\operatorname{id})\delta(w_1)
=g^{a'}\otimes w_1.
\]
Thus $\pi(c_{\mathbf r})=0$.
Applying $\pi\otimes\operatorname{id}$
to \eqref{eq:A-c-eq} gives
\[
(\pi\otimes\operatorname{id})
\Delta(c_{\mathbf r})
=g^{a'}\otimes c_{\mathbf r},
\]
so every PBW term of $c_{\mathbf r}$
has terminal group factor $g^{a'}$.

For $\ell\ge1$, the coproduct of
$z^\ell y^{\mathbf s}g^{a'}$ contains
\[
\ell z^{\ell-1}y^{\mathbf s}g^{a'}
\otimes zg^{a'+t\sum_i s_im_i}.
\]
Since the second factors in the sum
of \eqref{eq:A-c-eq} have zero
positive $z$-degree, comparison of
PBW coefficients forces
$\ell=1$ and $\mathbf s=\mathbf0$.
Consequently,
\[
c_{\mathbf r}
=\gamma_{\mathbf r}zg^{a'}
+\sum_{\mathbf0<\mathbf s<\mathbf e}
\gamma_{\mathbf r,\mathbf s}
y^{\mathbf s}g^{a'}.
\]

Comparing the $zg^{a'}$-coefficients
in \eqref{eq:A-c-eq}, we obtain
\[
\delta\left(
\sum_{\mathbf r<\mathbf d}
\gamma_{\mathbf r}y^{\mathbf r}w_0
\right)
=g^{a'}\otimes
\sum_{\mathbf r<\mathbf d}
\gamma_{\mathbf r}y^{\mathbf r}w_0.
\]
Since the only nonzero
coaction-homogeneous vectors of
$V_0(a,b)$ are scalar multiples
of $w_0$, we have
$\gamma_{\mathbf r}=0$ for
$\mathbf r\ne\mathbf0$.
Writing $\gamma=\gamma_{\mathbf0}$,
the coaction degree, $g$-weight,
and \eqref{eq:di} give
\eqref{eq:A-gamma} whenever
$\gamma\ne0$.

We now eliminate the pure-$y$ terms.
Applying $\operatorname{id}\otimes\pi$
to \eqref{eq:A-c-eq} and using the
$g$-weight condition, a nonzero
$\gamma_{\mathbf r,\mathbf s}$
must satisfy
\begin{equation}\label{eq:A-pure-cong}
\begin{aligned}
a-a'&\equiv
t\sum_i(r_i+s_i)m_i\pmod n,\\
b'-b&\equiv
-\frac{\sum_i(r_i+s_i)m_i}{m_0}
\pmod n.
\end{aligned}
\end{equation}
By the fraction bijection, there is
a unique reduced index $\mathbf k$
such that
\[
k_i\equiv r_i+s_i\pmod{e_i},
\qquad 0\le k_i<e_i.
\]
Thus, for each $\mathbf r$, at most
one reduced $\mathbf s$ can occur.

Comparing the coproduct coefficients
in \eqref{eq:A-c-eq} yields
\begin{equation}\label{eq:A-rec}
\gamma_{\mathbf r,\mathbf s}
[s_i]_{\xi^{-tm_i^2/m_0}}
=
\gamma_{\mathbf r+\mathbf1_i,
\mathbf s-\mathbf1_i}
C_{\mathbf r+\mathbf1_i,\mathbf1_i}(a,b),
\end{equation}
whenever $|\mathbf s|\ge2$,
$s_i>0$, and $r_i+1<d_i$.
The Gaussian integer is nonzero
for $0<s_i<e_i$.

First suppose that a linear term
occurs at an interior index
$\mathbf r=\mathbf k-\mathbf1_i$,
where $\mathbf k<\mathbf d$.
Write its coefficient as $\gamma_i$.
Comparison of the mixed
$y_i,y_j$-terms gives
\[
\gamma_i C_{\mathbf k,\mathbf1_j}(a,b)
=
\gamma_j C_{\mathbf k,\mathbf1_i}(a,b).
\]
All the indicated coefficients
$C_{\mathbf k,\mathbf1_i}(a,b)$
are nonzero. Hence there is
$c\in\mathbbm{k}$ such that
\[
\gamma_i=-cC_{\mathbf k,\mathbf1_i}(a,b).
\]
Replacing $w_1$ by
$w_1+cy^{\mathbf k}w_0$
eliminates all interior linear
terms simultaneously. The new
lift retains its $g$-weight and
projected coaction degree by
\eqref{eq:A-pure-cong}.

Next suppose a nonzero linear
term remains at the truncation
boundary
\[
k_i=d_i<e_i,\qquad
k_j<d_j\quad(j\ne i).
\]
Repeated application of
\eqref{eq:A-rec} gives
$\gamma_{\mathbf0,\mathbf k}\ne0$,
since all coefficients encountered
before the first zero are nonzero.
The congruences
\eqref{eq:A-pure-cong} and
\eqref{eq:di} imply
\[
d_i'=e_i-d_i,\qquad
1\le d_i'<e_i.
\]

In \eqref{eq:A-master}, compare
the coefficient of
\[
z\prod_{j\ne i}y_j^{k_j}
g^{a'-td_i'm_i}\otimes w_0.
\]
The right-hand side has zero
coefficient. On the left, only
$c_{\mathbf0}\otimes w_0$ can
contribute, because the action
of $\mathcal D_i$ never decreases
the $y_i$-degree of the second
tensor factor. Iterating
\eqref{eq:D}, and using
$k_i+d_i'=e_i$, gives
\[
0=
\gamma_{\mathbf0,\mathbf k}
\xi^{btm_id_i'}
\prod_{j=1}^{d_i'}
\left(1-\xi^{jtm_i^2/m_0}\right).
\]
The right-hand side is nonzero:
$\gamma_{\mathbf0,\mathbf k}\ne0$
and $\xi^{tm_i^2/m_0}$ has order
$e_i>d_i'$. This is a
contradiction.

It remains to consider the
periodic boundary
\[
k_i=0,\qquad d_i=e_i>1,
\qquad k_j<d_j\quad(j\ne i).
\]
A nonzero linear term at
$\mathbf k+(e_i-1)\mathbf1_i$
would, by \eqref{eq:A-rec},
force
\[
\gamma_{\mathbf k+\mathbf1_i,
(e_i-1)\mathbf1_i}\ne0,
\]
since all intermediate Gaussian
integers and first-zero
coefficients are nonzero.

Now consider \eqref{eq:A-c-eq}
with $\mathbf r=\mathbf k$.
The left-hand side has zero
coefficient at
\[
y_i^{e_i-1}g^{a'}\otimes
y_i g^{a-t(\sum_jk_jm_j+m_i)},
\]
because every reduced pure-$y$
term of $c_{\mathbf k}$ has
$i$th exponent zero.
On the right, the term indexed
by $\mathbf k+\mathbf1_i$
contributes
\[
\gamma_{\mathbf k+\mathbf1_i,
(e_i-1)\mathbf1_i}
C^{(i)}_{1,1}(a,b)\ne0.
\]
This contradiction excludes
the periodic-boundary terms.

Thus every linear pure-$y$
coefficient has been eliminated.
If a higher term
$\gamma_{\mathbf r,\mathbf s}$
with $|\mathbf s|\ge2$ remained,
comparison of the first factor
$y_i g^{a'}$, for some $s_i>0$,
in \eqref{eq:A-c-eq} would give
the nonzero coefficient
\[
\gamma_{\mathbf r,\mathbf s}
[s_i]_{\xi^{-tm_i^2/m_0}}
\]
on the left, whereas the right
would require a linear pure-$y$
coefficient. Hence all higher
terms vanish as well, proving
the coaction formula
in \eqref{eq:A-normal}.

Finally, $y_i^{d_i'}w_1$
belongs to $V_0(a,b)$.
The first-zero property gives
\[
\delta(y_i^{d_i'}w_1)
=g^{a'-td_i'm_i}
\otimes y_i^{d_i'}w_1
\]
when $\gamma=0$.
If $\gamma\ne0$, then
$a'=a$, $d_i'=d_i$, and
\[
\mathcal D_i^{d_i}(zg^a\otimes w_0)
=(z\otimes\operatorname{id})
\delta(y_i^{d_i}w_0)=0,
\]
so the same identity holds.
Since the only nonzero
coaction-homogeneous vectors
of $V_0(a,b)$ are scalar
multiples of $w_0$, we obtain
\[
y_i^{d_i'}w_1=u_iw_0.
\]
If $u_i\ne0$, comparison of the
coaction degrees and $g$-weights
gives \eqref{eq:A-u}.
\end{proof}

\begin{theorem}\label{thm:A-class}
Assume $m>1$ and consider the extension
\eqref{eq:A-extension}. By Lemma~\ref{lem:A-normal},
we may choose $w_1$ such that
\[
\delta(w_1)=g^{a'}\otimes w_1
+\gamma zg^a\otimes w_0,
\qquad
y_i^{d_i'}w_1=u_iw_0
\quad(1\le i\le\theta),
\]
where $u_i,\gamma\in\mathbbm{k}$.

The extension splits if and only if
\begin{equation}\label{eq:A-thm-split}
u_i=0\quad(1\le i\le\theta),
\qquad \gamma=0.
\end{equation}

For a nonsplit extension, the following
conditions hold.

\begin{enumerate}[
label=\textup{(\arabic*)},leftmargin=2.5em]

\item Suppose $u_i\ne0$ for some $i$.
Then
\begin{equation}\label{eq:A-thm-u}
\begin{aligned}
a&\equiv a'-td_i'm_i\pmod n,\\
b&\equiv b'-\frac{d_i'm_i}{m_0}
\pmod n.
\end{aligned}
\end{equation}
Moreover, we have the following two cases.

\begin{enumerate}[label=\textup{(\alph*)}]
\item If $d_i'<e_i$, then
\begin{equation}\label{eq:A-thm-proper}
\begin{aligned}
d_i&=e_i-d_i',&
d_j&=d_j'\quad(j\ne i),\\
u_j&=0\quad(j\ne i),&
\gamma&=0.
\end{aligned}
\end{equation}

\item If $d_i'=e_i$, then
\begin{equation}\label{eq:A-thm-full}
\begin{gathered}
d_j=d_j'=e_j,\qquad
u_j=u\ne0\quad(1\le j\le\theta),\\
a'=a,\qquad
b'-b\equiv\frac{e_1m_1}{m_0}\pmod n.
\end{gathered}
\end{equation}
\end{enumerate}

\item Suppose $\gamma\ne0$. Then
\begin{equation}\label{eq:A-thm-gamma}
\begin{gathered}
a'=a,\qquad
b'-b\equiv-\frac{e_1m_1}{m_0}\pmod n,\\
d_i'=d_i\qquad(1\le i\le\theta).
\end{gathered}
\end{equation}
If $d_i<e_i$ for some $i$, then
\begin{equation}\label{eq:A-thm-nonfull}
u_j=0\qquad(1\le j\le\theta).
\end{equation}
If $d_i=e_i$ for every $i$, then
$u_j=u$ for all $j$, and
\begin{equation}\label{eq:A-thm-mixed}
u\ne0\Longrightarrow
\frac{2e_1m_1}{m_0}\equiv0\pmod n.
\end{equation}
In particular, conditions
\textup{(1)(b)} and \textup{(2)}
may hold simultaneously, giving mixed
extensions with $u\gamma\ne0$.
\end{enumerate}

Conversely, every choice of parameters
satisfying the above conditions is
realizable by a Yetter-Drinfeld extension.

The corresponding middle modules have
basis \eqref{eq:A-basis}, actions
\eqref{eq:A-g-action}--\eqref{eq:A-middle-action},
and coactions specified below.
\end{theorem}

\begin{proof}
We first prove necessity. Suppose $u_i\ne0$.
By \eqref{eq:A-thm-u} and the first-zero
congruences \eqref{eq:di},
\begin{equation}\label{eq:A-d-compare}
m_j(d_j-d_j')\equiv-2d_i'm_i\pmod{e_j}.
\end{equation}
Since $e_j\mid m_i$ for $j\ne i$ and
$\gcd(e_j,m_j)=1$, we obtain
$d_j=d_j'$ for $j\ne i$. For $j=i$,
the same congruence gives
$d_i+d_i'\equiv0\pmod{e_i}$.

If $d_i'<e_i$, then $d_i=e_i-d_i'$.
If some $u_j\ne0$ with $j\ne i$,
comparison of \eqref{eq:A-thm-u}
for $i$ and $j$ gives
$m\mid d_i'm_i-d_j'm_j$, hence
$e_i\mid d_i'$, a contradiction.
Moreover, \eqref{eq:A-thm-u} implies
$a\not\equiv a'\pmod n$, whereas
$\gamma\ne0$ would imply $a'=a$
by \eqref{eq:A-gamma}. Thus
$\gamma=0$, proving
\eqref{eq:A-thm-proper}.

If $d_i'=e_i$, the common-power
relations give
\[
u_iw_0=zw_1
=u_jy_j^{e_j-d_j'}w_0
\qquad(1\le j\le\theta).
\]
Since $u_i\ne0$, the truncated
PBW basis forces $d_j'=e_j$
and $u_j=u_i$ for every $j$.
Equation \eqref{eq:A-d-compare}
then gives $d_j=e_j$, while
\eqref{eq:A-thm-u} and the fraction
conditions yield the remaining
relations in \eqref{eq:A-thm-full}.

Now suppose $\gamma\ne0$.
Lemma~\ref{lem:A-normal} gives
\eqref{eq:A-thm-gamma}.
If $u_j\ne0$, then
\eqref{eq:A-thm-u} implies
$e_j\mid d_j'$, so $d_j'=e_j$.
If some $d_i<e_i$, we have
$d_i'=d_i$ and $u_i=0$. Hence
$zw_1=y_i^{e_i-d_i}y_i^{d_i}w_1=0$,
contradicting $zw_1=u_jw_0\ne0$.
This proves \eqref{eq:A-thm-nonfull}.

If all $d_i=e_i$, the common-power
relations give $u_i=u$ for every
$i$. When $u\ne0$, comparison of
\eqref{eq:A-thm-full} and
\eqref{eq:A-thm-gamma} yields
\eqref{eq:A-thm-mixed}.
This completes necessity.

We next prove sufficiency.
If all parameters vanish, take
the direct sum of the two endpoint
Yetter-Drinfeld modules.

Suppose $u_i\ne0$ and $d_i'<e_i$.
Let $v$ be the standard generator
of $M(a',b')$ and define
\[
E=M(a',b')\Big/
\left(T(zv)+\sum_{j\ne i}T(y_j^{d_j'}v)\right),
\qquad
w_1=\bar v,\quad
w_0=u_i^{-1}y_i^{d_i'}w_1.
\]
The denominator is a Yetter-Drinfeld
submodule by Lemma~\ref{lem:di}.
Using \eqref{eq:A-thm-u} and
\eqref{eq:A-thm-proper}, we obtain
$\delta(w_0)=g^a\otimes w_0$,
$gw_0=\xi^bw_0$, and
$y_j^{d_j}w_0=0$ for every $j$.
The reduced PBW basis gives
\[
Tw_0\cong V_0(a,b),\qquad
E/Tw_0\cong V_0(a',b'),
\]
so the required extension exists.

Suppose next that $u_i=u\ne0$
for all $i$, with
\eqref{eq:A-thm-full}. Take
\[
E=M(a',b')/T(z^2v),\qquad
w_1=\bar v,\quad w_0=u^{-1}zw_1.
\]
Since $z^2v$ is standard,
this is a Yetter-Drinfeld quotient
with the induced coaction.
Its PBW basis and the relations
$zw_1=uw_0$, $zw_0=0$ show that
its submodule $Tw_0$ and quotient
are the prescribed simple modules.
This realizes the case $\gamma=0$.

It remains to consider $\gamma\ne0$.
If all $u_i=0$, take the direct sum
of the underlying left $T$-modules.
If $u_i=u\ne0$, use instead the
full-truncation underlying $T$-module
constructed above. In either case,
prescribe
\[
\delta(w_0)=g^a\otimes w_0,\qquad
\delta(w_1)=g^a\otimes w_1+
\gamma zg^a\otimes w_0.
\]
Coassociativity follows from
\[
\Delta(zg^a)
=g^a\otimes zg^a+zg^a\otimes g^a,
\]
and the $g$-compatibility follows
from \eqref{eq:A-thm-gamma}.
Moreover,
\[
\mathcal D_i^{d_i}(zg^a\otimes w_0)
=(z\otimes\operatorname{id})
\delta(y_i^{d_i}w_0)=0,
\]
so the truncation relations are
preserved. When $u\ne0$, we also have
$\delta(zw_1)=\delta(uw_0)$, since
$zw_0=0$. Thus the prescribed
coaction extends to a well-defined
Yetter-Drinfeld structure, realizing
all admissible parameters.

Finally, the zero parameters clearly
give a split extension. Conversely,
if the extension splits, there is a
standard lift $w_1+s$ with
$s\in V_0(a,b)$ satisfying the
quotient's truncation relations.
Since $\delta(s)$ has no positive
$z$-degree in its first PBW factor,
the additional term
$\gamma zg^a\otimes w_0$ cannot be
cancelled. Hence $\gamma=0$.
The standard-element property then
implies $s=0$ or $s=cw_0$ with
$(a',b')=(a,b)$. In either case,
$y_i^{d_i'}s=0$, so
$0=y_i^{d_i'}(w_1+s)=u_iw_0$
for every $i$. Thus all $u_i=0$,
proving \eqref{eq:A-thm-split}.
\end{proof}

For completeness, the $T$-action on the basis \eqref{eq:A-basis}
is
\begin{equation}\label{eq:A-g-action}
 g(y^{\mathbf r}w_0)=\xi^{b-\sum_i r_im_i/m_0}y^{\mathbf r}w_0,
 \qquad g(y^{\mathbf s}w_1)=
 \xi^{b'-\sum_i s_im_i/m_0}y^{\mathbf s}w_1,
\end{equation}
\begin{equation}\label{eq:A-middle-action}
\begin{aligned}
 y_i(y^{\mathbf r}w_0)&=\begin{cases}
 y^{\mathbf r+\mathbf1_i}w_0,&r_i<d_i-1,\\
 0,&r_i=d_i-1,
 \end{cases}\\
 y_i(y^{\mathbf s}w_1)&=\begin{cases}
 y^{\mathbf s+\mathbf1_i}w_1,&s_i<d_i'-1,\\
 u_i\prod_{j\ne i}y_j^{s_j}w_0,&s_i=d_i'-1.
 \end{cases}
\end{aligned}
\end{equation}
An out-of-range lower-family monomial is interpreted as zero.
The coaction on the lower family is \eqref{eq:C-coaction}
with $w_0$ in place of $v$. If $\gamma=0$, the upper-family
coaction is the same formula with $(a',b',w_1)$ in place of
$(a,b,v)$. If $\gamma\ne0$, then
$\xi^{b'tm_i}=\xi^{btm_i}$ and hence
$C_{\mathbf r,\mathbf j}(a,b')=C_{\mathbf r,\mathbf j}(a,b)$;
therefore
\begin{equation}\label{eq:A-full-coaction}
\begin{aligned}
 \delta(y^{\mathbf r}w_1)
 ={}&\sum_{\mathbf0\le\mathbf j\le\mathbf r}
 C_{\mathbf r,\mathbf j}(a,b)y^{\mathbf j}
 g^{a-t\sum_i r_im_i}\otimes y^{\mathbf r-\mathbf j}w_1\\
 &+\gamma\sum_{\mathbf0\le\mathbf j\le\mathbf r}
 C_{\mathbf r,\mathbf j}(a,b)zy^{\mathbf j}
 g^{a-t\sum_i r_im_i}\otimes y^{\mathbf r-\mathbf j}w_0.
\end{aligned}
\end{equation}
In case (1) with $u_i\ne0$, $E$ may also be described using
its single-generator basis displayed in the proof. In case (2)
with $u\ne0$, the single-generator basis is
$\{z^\ell y^{\mathbf r}w_1:\ell=0,1,\ 0\le r_i<e_i\}$.

\begin{proposition}\label{prop:A-iso-ext}
For fixed ordered endpoints $V_0(a,b)$
and $V_0(a',b')$, the admissible
parameters $u_i$ ($1\le i\le\theta$)
and $\gamma$ in Theorem~\ref{thm:A-class}
give a linear parametrization of
$\operatorname{Ext}^1(V_0(a',b'),V_0(a,b))$.

Every nonsplit middle module is
indecomposable. Moreover, for fixed
ordered endpoints, two nonsplit middle
modules are isomorphic if and only if
\[
E(u_i,\gamma)\cong E(u_i',\gamma')
\iff
u_i'=cu_i\ (1\le i\le\theta),
\quad \gamma'=c\gamma
\]
for some $c\in\mathbbm{k}^{\times}$.

When both the full-truncation action
and coaction directions are admissible,
their nonsplit middle modules are
parametrized by
$[u:\gamma]\in\mathbb P^1(\mathbbm{k})$.
\end{proposition}

\begin{proof}
We first show that the parameters are
independent of the normalized lift.
Let $w_1^*=w_1+s$, where
$s\in V_0(a,b)$, be another lift
in the normal form of Lemma~\ref{lem:A-normal}.
Comparing their coactions gives
\[
\delta(s)-g^{a'}\otimes s
=(\gamma^*-\gamma)zg^a\otimes w_0.
\]
Since the first coaction factors of
$V_0(a,b)$ have zero positive
$z$-degree, we obtain
$\gamma^*=\gamma$ and
$\delta(s)=g^{a'}\otimes s$.
Moreover, $gs=\xi^{b'}s$.

The standard-element property of
$V_0(a,b)$ implies that either
$s=0$, or $s=cw_0$ with
$(a',b')=(a,b)$. In the latter case,
$d_i'=d_i$ and $y_i^{d_i'}s=0$.
Thus $u_i^*=u_i$ for every $i$.

For the Baer sum, choose the pair
of normalized upper lifts in the
pullback and apply the pushout
along the addition map on the
lower simple modules. The boundary
actions and added coaction terms
then add componentwise.
Scalar multiplication acts similarly.
Together with Theorem~\ref{thm:A-class},
this identifies the extension classes
with the admissible parameter space.

By Theorem~\ref{thm:A-class}, a nonzero admissible parameter
choice gives a nonsplit extension. Lemma~\ref{lem:extension-indecomp}
then shows that its middle module is indecomposable and that
the lower simple is its unique nonzero proper submodule.

Now let
$f:E(u_i,\gamma)\to E(u_i',\gamma')$
be an isomorphism between nonsplit
middle modules. The uniqueness of
the lower simple and Schur's lemma give
\[
f(w_0)=\alpha w_0',\qquad
f(w_1)=\beta w_1'+s,
\]
where $\alpha,\beta\in\mathbbm{k}^\times$
and $s$ belongs to the lower simple.

Colinearity gives
\[
\delta(s)-g^{a'}\otimes s
=(\alpha\gamma-\beta\gamma')zg^a\otimes w'_0.
\]
Since the first coaction factors of the lower simple have
zero positive $z$-degree, we obtain
\[
\alpha\gamma=\beta\gamma',
\qquad
\delta(s)=g^{a'}\otimes s.
\]
Together with $gs=\xi^{b'}s$, the standard-element property
implies that $s=0$, or $s=cw'_0$ with
$(a',b')=(a,b)$. Hence $y_i^{d_i'}s=0$ for every $i$.
Comparing the boundary actions gives
\[
\alpha u_i=\beta u_i'
\qquad (1\le i\le\theta).
\]
Thus the parameters are proportional.

Conversely, if these equalities hold,
the assignments
$f(w_0)=\alpha w_0'$ and
$f(w_1)=\beta w_1'$ respect the
actions and coactions and extend
to a Yetter-Drinfeld isomorphism.
The projective parametrization
follows immediately.
\end{proof}

\begin{corollary}\label{cor:A-ext}
For $m>1$,
\begin{equation}\label{eq:A-ext-dim}
\begin{aligned}
&\dim_{\mathbbm{k}}\operatorname{Ext}^1(V_0(a',b'),V_0(a,b))\\
={}&\sum_{\{i:d_i'<e_i\}}
 \mathbf1_{\{a\equiv a'-td_i'm_i,\ 
 b\equiv b'-d_i'm_i/m_0\}}\\
&+\mathbf1_{\{d_j'=e_j\ \forall j,\ a'=a,\ 
 b'-b\equiv e_1m_1/m_0\}}\\
&+\mathbf1_{\{a'=a,\ b'-b\equiv-e_1m_1/m_0\}}.
\end{aligned}
\end{equation}
Here $\mathbf1_{\{\cdots\}}$ denotes the indicator of a stated
condition. In particular this dimension is at most $2$.
\end{corollary}

\begin{proof}
By Theorem~\ref{thm:A-class} and
Proposition~\ref{prop:A-iso-ext},
the three terms in \eqref{eq:A-ext-dim}
count the admissible independent
directions arising from proper
truncation, full-truncation action,
and nonzero coaction, respectively.

The first-zero congruences
\eqref{eq:di} show that the endpoint
conditions in the last two terms
already imply the required relations
among $d_i$ and $d_i'$.

At most one index contributes to the
first sum. Indeed, if two distinct
indices $i,j$ contributed, then
$m\mid d_i'm_i-d_j'm_j$, which,
using $e_i\mid m_j$ and
$\gcd(e_i,m_i)=1$, would imply
$e_i\mid d_i'$, a contradiction.

The proper-truncation direction
cannot coexist with either of the
other directions. Hence at most
the full-truncation action and
coaction directions occur together,
and the dimension is at most $2$.
\end{proof}

\subsection{Extensions between cyclic simple modules}

Consider
\begin{equation}\label{eq:B-extension}
 0\longrightarrow V_\lambda(a,b)=Tw_0\longrightarrow E
 \xrightarrow{p}V_\mu(a',b')=Tp(w_1)\longrightarrow0,
 \qquad\lambda,\mu\in\mathbbm{k}^\times.
\end{equation}
By Corollary~\ref{cor:unequal}, it splits when $\lambda\ne\mu$.
Suppose henceforth $\lambda=\mu$. As $z^t-\lambda$ vanishes
on both endpoints, it induces a Yetter-Drinfeld morphism
\begin{equation}\label{eq:connecting}
 V_\lambda(a',b')\longrightarrow V_\lambda(a,b),\qquad
 p(v)\longmapsto(z^t-\lambda)v.
\end{equation}
Thus $(z^t-\lambda)E=0$ if the endpoints are nonisomorphic.
For isomorphic endpoints, identify both with the same model
$S=V_\lambda(a,b)$; the distinct $g$-weights of $S$ then give
\begin{equation}\label{eq:B-action}
 gw_1=\xi^bw_1,\qquad(z^t-\lambda)w_1=uw_0
 \quad(u\in\mathbbm{k}).
\end{equation}

Choose a lift $w_1$ in the $g^{a'}$-component of the projected
coaction and in the $\xi^{b'}$-eigenspace of $g$. It has coaction
\begin{equation}\label{eq:B-c}
 \delta(w_1)=g^{a'}\otimes w_1+
 \sum_{\substack{0\le\ell<t\\0\le r_i<e_i}}
 c_{\ell,\mathbf r}\otimes z^\ell y^{\mathbf r}w_0.
\end{equation}
The coefficient equations are
\begin{equation}\label{eq:B-c-eq}
\begin{aligned}
 \varepsilon(c_{\ell,\mathbf r})&=0,\\
 \Delta(c_{\ell,\mathbf r})&=g^{a'}\otimes c_{\ell,\mathbf r}
 +\sum_{\mathbf r\le\mathbf s<\mathbf e}
 C_{\mathbf s,\mathbf s-\mathbf r}(a,b)c_{\ell,\mathbf s}
 \otimes y^{\mathbf s-\mathbf r}g^{a-t\sum_i s_im_i},\\
 gc_{\ell,\mathbf r}g^{-1}&=
 \xi^{b'-b+(\ell e_1m_1+\sum_i r_im_i)/m_0}
 c_{\ell,\mathbf r}.
\end{aligned}
\end{equation}
The index $\ell$ is unchanged in the second line because
$\delta(z^\ell v)=(\operatorname{id}\otimes z^\ell)\delta(v)$ on
$V_\lambda(a,b)$.

\begin{lemma}\label{lem:B-reduced}
Every coefficient $c_{\ell,\mathbf r}$
in \eqref{eq:B-c} is a linear
combination of $zg^{a'}$ and
$y^{\mathbf s}g^{a'}$, where
$0\le s_i<e_i$ and
$\mathbf s\ne\mathbf0$.

A nonzero pure-$y$ coefficient
requires
\begin{equation}\label{eq:B-cong}
\begin{aligned}
a-a'&\equiv
t\sum_i(r_i+s_i)m_i\pmod n,\\
b'-b&\equiv
-\frac{\ell e_1m_1+
\sum_i(r_i+s_i)m_i}{m_0}
\pmod n.
\end{aligned}
\end{equation}
A nonzero $zg^{a'}$ coefficient
requires
\begin{equation}\label{eq:B-z-cong}
\begin{aligned}
a-a'&\equiv t\sum_i r_im_i\pmod n,\\
b'-b&\equiv
-\frac{(\ell+1)e_1m_1+
\sum_i r_im_i}{m_0}\pmod n.
\end{aligned}
\end{equation}
\end{lemma}

\begin{proof}
Let $\pi:T\to\mathbbm{k}\langle g\rangle$
be the Hopf projection with
$\pi(y_i)=0$. By the normalization
of $w_1$, we have
$\pi(c_{\ell,\mathbf r})=0$.
Applying the two projections to
\eqref{eq:B-c-eq} gives
\[
\begin{aligned}
(\pi\otimes\operatorname{id})
\Delta(c_{\ell,\mathbf r})
&=g^{a'}\otimes c_{\ell,\mathbf r},\\
(\operatorname{id}\otimes\pi)
\Delta(c_{\ell,\mathbf r})
&=c_{\ell,\mathbf r}\otimes
g^{a-t\sum_i r_im_i}.
\end{aligned}
\]
The first identity implies that
every PBW term in $c_{\ell,\mathbf r}$
has terminal group factor $g^{a'}$.

For $k\ge1$, the coproduct of
$z^k y^{\mathbf s}g^{a'}$
contains
\[
kz^{k-1}y^{\mathbf s}g^{a'}
\otimes zg^{a'+t\sum_i s_im_i}.
\]
Since the summation in
\eqref{eq:B-c-eq} has no positive
$z$-degree in its second tensor
factors, PBW comparison forces
$k=1$ and $\mathbf s=\mathbf0$.
Together with $\pi(c_{\ell,\mathbf r})=0$,
this yields
\[
c_{\ell,\mathbf r}
=\gamma_{\ell,\mathbf r}zg^{a'}
+\sum_{\substack{
0\le s_i<e_i\\\mathbf s\ne\mathbf0}}
\gamma_{\ell,\mathbf r,\mathbf s}
y^{\mathbf s}g^{a'}.
\]

Finally,
\[
\begin{aligned}
(\operatorname{id}\otimes\pi)
\Delta(y^{\mathbf s}g^{a'})
&=y^{\mathbf s}g^{a'}\otimes
g^{a'+t\sum_i s_im_i},\\
(\operatorname{id}\otimes\pi)
\Delta(zg^{a'})
&=zg^{a'}\otimes g^{a'}.
\end{aligned}
\]
Comparing these identities with
the projected coproduct and
$g$-weight equations in
\eqref{eq:B-c-eq} gives
\eqref{eq:B-cong} and
\eqref{eq:B-z-cong}.
\end{proof}

\begin{proposition}\label{prop:B-noniso}
If $V_\lambda(a,b)\not\cong V_\lambda(a',b')$, then every
extension \eqref{eq:B-extension} splits.
\end{proposition}

\begin{proof}
Since $z^t-\lambda$ annihilates
both endpoints, it induces the
Yetter-Drinfeld morphism
\eqref{eq:connecting}. By Schur's
lemma and the nonisomorphism
assumption, this morphism is zero.
Hence
\begin{equation}\label{eq:B-noniso-central}
(z^t-\lambda)E=0.
\end{equation}

Write the coefficients in
\eqref{eq:B-c} as in
Lemma~\ref{lem:B-reduced}.
Suppose some $zg^{a'}$ coefficient
is nonzero, and set
\[
v=\sum_{\ell,\mathbf r}
\gamma_{\ell,\mathbf r}
z^\ell y^{\mathbf r}w_0\ne0.
\]
Comparing the $zg^{a'}$-terms
in \eqref{eq:B-c-eq} gives
$\delta(v)=g^{a'}\otimes v$,
while the $g$-weight equation yields
$gv=\xi^{b'+e_1m_1/m_0}v$.

Since $z^t=\lambda\ne0$ on the
lower cyclic simple, $zv\ne0$.
Moreover,
\[
\delta(zv)=g^{a'}\otimes zv,
\qquad g(zv)=\xi^{b'}zv.
\]
Thus $zv$ is standard of type
$(a',b')$, contradicting
Proposition~\ref{prop:iso}.
Consequently,
\begin{equation}\label{eq:B-no-z}
\gamma_{\ell,\mathbf r}=0
\qquad\text{for all }\ell,\mathbf r.
\end{equation}

We next consider the linear
pure-$y$ coefficients.
By \eqref{eq:B-cong}, any nonzero
$y_i g^{a'}$ term determines a
nonzero reduced PBW vector
$v=z^Ly^{\mathbf R}w_0$ satisfying
\[
gv=\xi^{b'}v,\qquad
(\pi\otimes\operatorname{id})\delta(v)
=g^{a'}\otimes v.
\]
Such a vector is unique, since
the reduced PBW vectors have
distinct $g$-weights.
Thus all linear terms have the form
\begin{equation}\label{eq:B-linear-v}
\sum_{\{i:e_i>1\}}
\gamma_i y_i g^{a'}
\otimes y_i^{-1}v,
\qquad
y_i^{-1}=\lambda^{-1}z^{t-1}y_i^{e_i-1}.
\end{equation}
If no such vector exists, all
linear coefficients vanish.

We claim that
\begin{equation}\label{eq:B-cross-rel}
\gamma_iC^{(j)}_{R_j,1}(a,b)
=\gamma_jC^{(i)}_{R_i,1}(a,b)
\qquad(i\ne j,\ e_i,e_j>1).
\end{equation}
If $R_i,R_j>0$, this follows by
comparing the two mixed terms
\[
\begin{aligned}
&y_i g^{a'}\otimes
y_jg^{a'+tm_i},\\
&y_j g^{a'}\otimes
y_ig^{a'+tm_j}
\end{aligned}
\]
in the coproduct of
$y_iy_jg^{a'}$ and in
\eqref{eq:B-c-eq}. Their
coproduct coefficients are both
$1$ by the fraction condition.

The same comparison applies at
the periodic boundary.
Suppose $R_i=0$ and $R_j>0$.
After reducing $y_i^{-1}v$ by
$y_i^{e_i}=z$ and $z^t=\lambda$,
its $i$th exponent is $e_i-1$,
with a nonzero scalar factor
$1$ or $\lambda^{-1}$.
Consider \eqref{eq:B-c-eq} at
the corresponding reduced index
whose $j$th exponent is $R_j-1$.

Comparing the mixed PBW coefficients in
\eqref{eq:B-c-eq} gives
\[
\kappa\gamma_i C^{(j)}_{R_j,1}(a,b)=0,
\]
where $\kappa\in\mathbbm{k}^\times$ arises from
the reduction using $z^t=\lambda$.
Indeed, the opposite mixed term cannot occur
on the right, since it would require a reduced
index with $i$th exponent $e_i$.
Hence
\[
\gamma_i C^{(j)}_{R_j,1}(a,b)
=0=\gamma_j C^{(i)}_{0,1}(a,b),
\]
which proves \eqref{eq:B-cross-rel}.

The case $R_j=0$ is symmetric;
if both exponents vanish, the
relation is immediate.

If every $C^{(i)}_{R_i,1}(a,b)$
vanishes, Lemma~\ref{lem:di}
implies $R_i\in\{0,d_i\}$
for every $i$, with $R_i=0$
when $d_i=e_i$.
Hence $v$ is standard of type
$(a',b')$, contradicting
Proposition~\ref{prop:iso}.

We may therefore choose $j$ with
$C^{(j)}_{R_j,1}(a,b)\ne0$.
By \eqref{eq:B-cross-rel},
there exists $c\in\mathbbm{k}$
such that
\[
\gamma_i=cC^{(i)}_{R_i,1}(a,b)
\qquad\text{for all }i.
\]
Replacing $w_1$ by $w_1-cv$
eliminates all linear pure-$y$
terms. It does not introduce
$z$-terms, by the known coaction
of the lower cyclic simple.

If a higher pure-$y$ coefficient
remained, its coproduct would
have a nonzero term with first
factor $y_i g^{a'}$ for some $i$.
Equation \eqref{eq:B-c-eq}
would then require a linear
coefficient, a contradiction.
Thus the new lift satisfies
\[
\delta(w_1)=g^{a'}\otimes w_1,
\qquad gw_1=\xi^{b'}w_1.
\]
Together with
\eqref{eq:B-noniso-central},
this gives
$(z^t-\lambda)w_1=0$.

Hence the induced morphism
$M(a',b')\to E$ factors through
$V_\lambda(a',b')$ and defines
a Yetter-Drinfeld section of
the quotient map. Therefore
the extension splits.
\end{proof}

\begin{proposition}\label{prop:B-normal}
Let $S=V_\lambda(a,b)$. Every self-extension admits a lift
satisfying \eqref{eq:B-action} and
\begin{equation}\label{eq:B-normal}
\begin{aligned}
 \delta(w_1)={}&g^a\otimes w_1
 +\gamma zg^a\otimes z^{t-1}w_0\\
 &+\sum_{\{i:d_i<e_i\}}\gamma_i
 \sum_{j=1}^{e_i-d_i}
 \frac{C^{(i)}_{e_i-1,j-1}(a,b)}
 {[j]_{\xi^{-tm_i^2/m_0}}}
 y_i^jg^a\otimes z^{t-1}y_i^{e_i-j}w_0,
\end{aligned}
\end{equation}
where $\gamma,\gamma_i\in\mathbbm{k}$ are independent parameters.
\end{proposition}

\begin{proof}
Taking $a'=a$ and $b'=b$ in
Lemma~\ref{lem:B-reduced}, the fraction
bijection and the root compatibility
give the following restrictions on a
nonzero pure-$y$ coefficient:
\begin{equation}\label{eq:B-self-cong}
r_i+s_i\in\{0,e_i\}\quad(1\le i\le\theta),
\qquad
\ell+\#\{i:s_i>0\}\equiv0\pmod t.
\end{equation}

If two distinct coordinates $s_i,s_j$
were positive, comparison of the tensor
\[
y_i g^a\otimes
y^{\mathbf s-\mathbf1_i}g^{a+tm_i}
\]
in \eqref{eq:B-c-eq} would require a
reduced index whose $j$th coordinate
is $r_j+s_j=e_j$, a contradiction.
Thus $\mathbf s$ has only one positive
coordinate.

The congruences for the $z$-coefficient
are treated similarly. Consequently,
the only possible additional coefficients
are
\begin{equation}\label{eq:B-self-candidates}
\begin{aligned}
c_{t-1,\mathbf0}&=\gamma zg^a,\\
c_{t-1,(e_i-s)\mathbf1_i}
&=\gamma_{i,s}y_i^sg^a,
\qquad 1\le s<e_i.
\end{aligned}
\end{equation}

Comparing the coefficient of
$y_i^{s-1}g^a\otimes
y_i g^{a+t(s-1)m_i}$
in \eqref{eq:B-c-eq}, we obtain
\begin{equation}\label{eq:B-rec}
\gamma_{i,s}[s]_{\xi^{-tm_i^2/m_0}}
=\gamma_{i,s-1}C^{(i)}_{e_i-s+1,1}(a,b),
\qquad 2\le s<e_i.
\end{equation}
Since the Gaussian integers in \eqref{eq:B-rec} are nonzero,
the product formula \eqref{eq:C} and the recurrences
\eqref{eq:Gaussian-recurrence} give
\begin{equation}\label{eq:B-self-solution}
\gamma_{i,s}
=\gamma_i
\frac{C^{(i)}_{e_i-1,s-1}(a,b)}
{[s]_{\xi^{-tm_i^2/m_0}}},
\qquad \gamma_i=\gamma_{i,1}.
\end{equation}

At the periodic boundary
$\mathbf r=\mathbf0$, comparison of
the coefficient of
$y_i^{e_i-1}g^a\otimes
y_i g^{a+t(e_i-1)m_i}$
gives
\[
\gamma_i
\prod_{j=1}^{e_i-1}
C^{(i)}_{j,1}(a,b)=0,
\]
up to nonzero Gaussian factors.
By Lemma~\ref{lem:di}, this forces
$\gamma_i=0$ when $d_i=e_i$.
When $d_i<e_i$, formula
\eqref{eq:B-self-solution} gives
$\gamma_{i,s}=0$ for $s>e_i-d_i$.

It remains to verify the full
coassociativity equations.
For $1\le j<s<e_i$, they reduce to
\[
\gamma_{i,s}
{s\brack j}_{\xi^{-tm_i^2/m_0}}
=
\gamma_{i,s-j}
C^{(i)}_{e_i-s+j,j}(a,b).
\]
Substituting \eqref{eq:B-self-solution} and
\eqref{eq:C}, this identity follows from
\eqref{eq:Gaussian-recurrence}.

The remaining boundary equations are
\[
\gamma_{i,e_i-j}C^{(i)}_{j,j}(a,b)=0,
\qquad 1\le j<e_i.
\]
If $j\ge d_i$, the second factor
vanishes; if $j<d_i$, the first
factor vanishes. The $z$-term is
primitive, and no mixed-direction
terms remain.

Thus all coassociativity equations
hold, and the freely chosen
coefficients are precisely
$\gamma$ and $\gamma_i$ for
$d_i<e_i$. This proves
\eqref{eq:B-normal}.
\end{proof}

\begin{theorem}\label{thm:B-class}
Let $S=V_\lambda(a,b)$. Up to extension equivalence, its
self-extensions are uniquely parametrized by $u,\gamma\in\mathbbm{k}$
and $\gamma_i\in\mathbbm{k}$ ($d_i<e_i$), through
\eqref{eq:B-action} and \eqref{eq:B-normal}. Every parameter
choice occurs, and the extension splits exactly when all
parameters vanish.

The middle module has dimension $2n$ and basis
\begin{equation}\label{eq:B-basis}
 \{z^\ell y^{\mathbf r}w_{0},\ z^\ell y^{\mathbf r}w_1:
 0\le\ell<t,\ 0\le r_i<e_i\}.
\end{equation}
On both families $g$ acts by \eqref{eq:g-on-basis}, while
\begin{equation}\label{eq:B-middle-action}
 y_i(z^\ell y^{\mathbf r}w_\varepsilon)=
 \begin{cases}
 z^\ell y^{\mathbf r+\mathbf1_i}w_\varepsilon,&r_i<e_i-1,\\
 z^{\ell+1}\prod_{j\ne i}y_j^{r_j}w_\varepsilon,
 &r_i=e_i-1,\ \ell<t-1,\\
 \lambda\prod_{j\ne i}y_j^{r_j}w_\varepsilon
 +\varepsilon u\prod_{j\ne i}y_j^{r_j}w_0,
 &r_i=e_i-1,\ \ell=t-1,
 \end{cases}
\end{equation}
for $\varepsilon\in\{0,1\}$. The lower coaction is given by
\eqref{eq:C-coaction}; the upper coaction is
\begin{equation}\label{eq:B-all-coaction}
\begin{aligned}
 \delta(z^\ell y^{\mathbf r}w_1)={}&
 \sum_{\mathbf0\le\mathbf j\le\mathbf r}
 C_{\mathbf r,\mathbf j}(a,b)y^{\mathbf j}
 g^{a-t\sum_i r_im_i}\otimes z^\ell y^{\mathbf r-\mathbf j}w_1\\
 &+(\operatorname{id}\otimes z^\ell)
 \mathcal D_1^{r_1}\cdots\mathcal D_\theta^{r_\theta}
 \left(\delta(w_1)-g^a\otimes w_1\right).
\end{aligned}
\end{equation}
Every nonzero parameter tuple yields an indecomposable middle
module. For fixed $S$, two such modules are isomorphic if and
only if their parameter tuples are proportional. Consequently,
\begin{equation}\label{eq:B-ext-dim}
 \dim_{\mathbbm{k}}\operatorname{Ext}^1(S,S)=
 2+\#\{i:d_i<e_i\}.
\end{equation}
\end{theorem}

\begin{proof}
By Proposition~\ref{prop:B-normal},
the additional coaction parameters
are precisely $\gamma$ and
$\gamma_i$ with $d_i<e_i$.
Since $z^t-\lambda$ induces an
endomorphism of the simple module
$S$, Schur's lemma gives the
action parameter $u$ in
\eqref{eq:B-action}.

We first construct the underlying
$T$-module. If $u=0$, take the
direct sum of two copies of $S$
as left $T$-modules. If $u\ne0$,
take
\[
E=M(a,b)/(z^t-\lambda)^2M(a,b),
\qquad
w_0=u^{-1}(z^t-\lambda)w_1.
\]
Freeness over $\mathbbm{k}[z^t]$
gives the basis \eqref{eq:B-basis}.
The algebra relations yield
\eqref{eq:B-middle-action}.

We prescribe the generator coactions
by \eqref{eq:B-normal}.
Proposition~\ref{prop:B-normal}
verifies their coassociativity and
counit conditions. The $g$-weight
condition follows from
$\xi^{te_im_i/m_0}=1$.

To extend the coaction, use the
adjoint-diagonal $T$-action on
$T\otimes E$, as in \eqref{eq:D}.
It remains to check that the
defining relations of the underlying
$T$-module are preserved.

Since $t\mid a$, every first tensor
factor in the additional coaction
\eqref{eq:B-normal} commutes with
$z$. As $z$ is primitive,
\[
\delta(z^tw_1)
=(\operatorname{id}\otimes z^t)
\delta(w_1).
\]
All additional second tensor factors
lie in $S$, where $z^t=\lambda$.
Consequently,
\[
\delta((z^t-\lambda)w_1)
=g^a\otimes(z^t-\lambda)w_1
=u g^a\otimes w_0.
\]
Thus the relation
$(z^t-\lambda)w_1=uw_0$
is compatible with the coaction.
The relation
$(z^t-\lambda)w_0=0$ follows
similarly, and the $g$-relations
follow from the weight condition.
Hence the prescribed coaction
descends to the underlying
$T$-module and satisfies the
Yetter-Drinfeld condition.
Formula \eqref{eq:B-all-coaction}
then follows by applying
\eqref{eq:D}.

We next consider extension
equivalence. Two normalized lifts
of the same quotient generator
differ by an element of $S$ having
$g$-weight $\xi^b$ and projected
coaction degree $g^a$.
Since the reduced PBW vectors of
$S$ have distinct $g$-weights,
this element is a scalar multiple
of $w_0$.
Such a change affects none of
$u,\gamma,\gamma_i$.

Baer addition acts componentwise
on these parameters. Together
with Proposition~\ref{prop:B-normal},
this gives
\eqref{eq:B-ext-dim}.
The extension splits exactly when
all parameters vanish.

Finally, Lemma~\ref{lem:extension-indecomp} shows that every
nonsplit middle module is indecomposable and has a unique
nonzero proper submodule.
As in Proposition~\ref{prop:A-iso-ext},
an isomorphism of nonsplit middle
modules acts by nonzero scalars
on the lower simple and the quotient.
Comparing \eqref{eq:B-action} and
\eqref{eq:B-normal} shows that
all parameters are multiplied by
one common nonzero scalar.
Conversely, such a rescaling
defines a Yetter-Drinfeld
isomorphism. This proves the
indecomposability and isomorphism
statements.
\end{proof}

\begin{corollary}\label{cor:B-ext}
For $\lambda,\mu\in\mathbbm{k}^\times$,
\begin{equation}\label{eq:B-all-ext-dim}
\dim_{\mathbbm{k}}\operatorname{Ext}^1(V_\mu(a',b'),V_\lambda(a,b))=
\begin{cases}
 2+\#\{i:d_i(a,b)<e_i\},&
 V_\mu(a',b')\cong V_\lambda(a,b),\\
 0,&\text{otherwise}.
\end{cases}
\end{equation}
In the nonzero case the explicit middle modules are obtained by
transport of the quotient identification to the same cyclic model.
\end{corollary}

\begin{proof}
If $\lambda\ne\mu$, the assertion
follows from Corollary~\ref{cor:unequal}.
For equal eigenvalues and
nonisomorphic cyclic simples,
it follows from
Proposition~\ref{prop:B-noniso}.

If the endpoints are isomorphic,
transporting the quotient along
an isomorphism identifies the
extension space with
$\operatorname{Ext}^1(S,S)$.
The dimension is therefore
given by \eqref{eq:B-ext-dim}.

To see that this number is independent
of the chosen presentation of $S$,
suppose
$V_\lambda(a,b)\cong V_\lambda(a',b')$.
By \eqref{eq:cyclic-iso} and the
first-zero congruences \eqref{eq:di},
\[
d_i(a',b')=
\begin{cases}
e_i-d_i(a,b),&i\in J,\\
d_i(a,b),&i\notin J,
\end{cases}
\]
for the subset $J$ appearing in
\eqref{eq:cyclic-iso}.
Thus the set of indices satisfying
$d_i<e_i$ is unchanged, proving
the assertion.
\end{proof}

\subsection{The case \texorpdfstring{$m=1$}{m=1} and the combined result}

If $m=1$, then $e_i=1$ and $y_i=z$ for every $i$;
moreover $t=n$. Every finite-dimensional truncated simple is
one-dimensional, of the form $V_0(0,b)=\mathbbm{k}w_0$, with
\[
 \delta(w_0)=1\otimes w_0,\qquad
 gw_0=\xi^bw_0,\qquad zw_0=0.
\]
Consider an extension
\[
0\longrightarrow V_0(0,b)\longrightarrow E
\longrightarrow V_0(0,b')\longrightarrow0.
\]
By Theorem~\ref{thm:primary}, the coaction of $E$
takes values in $\mathbbm{k}[z]\otimes E$.
Since $g$ is diagonalizable, we may choose
$gw_1=\xi^{b'}w_1$. Then
\[
 \delta(w_1)=1\otimes w_1+c\otimes w_0,
 \qquad zw_1=uw_0.
\]
Coassociativity and the counit imply
\[
 \Delta(c)=1\otimes c+c\otimes1,
 \qquad \varepsilon(c)=0.
\]
The PBW coproduct gives $P(T)=\mathbbm{k}z$,
so $c=\gamma z$. The $g$-action and
$g$-coaction compatibility yield
\begin{equation}\label{eq:m1}
\begin{aligned}
 zw_1&=uw_0,\quad zw_0=0,\qquad
 \delta(w_1)=1\otimes w_1+\gamma z\otimes w_0,\\
 u\ne0&\ \Longrightarrow\
 b'-b\equiv m_1/m_0\pmod n,\\
 \gamma\ne0&\ \Longrightarrow\
 b'-b\equiv-m_1/m_0\pmod n.
\end{aligned}
\end{equation}
Conversely, since $z$ is primitive and commutes
with $\mathbbm{k}[z]$, the Yetter-Drinfeld
compatibility gives
\[
 \delta(zw_1)
 =(\operatorname{id}\otimes z)\delta(w_1)
 =u(1\otimes w_0)
 =\delta(uw_0).
\]
Thus every pair $(u,\gamma)$ satisfying
\eqref{eq:m1} defines a Yetter-Drinfeld extension.
Changes of normalized lift preserve both parameters,
and the extension splits precisely when $u=\gamma=0$.
Consequently,
\[
\begin{aligned}
&\dim_{\mathbbm{k}}\operatorname{Ext}^1
 \bigl(V_0(0,b'),V_0(0,b)\bigr)\\
&\qquad=
\mathbf{1}_{\,b'-b\equiv m_1/m_0\pmod n}
+\mathbf{1}_{\,b'-b\equiv-m_1/m_0\pmod n},
\end{aligned}
\]
where $\mathbf{1}_P$ equals $1$ if $P$ holds
and $0$ otherwise.

Since $e_i=1$ for all $i$,
Theorem~\ref{thm:B-class} gives
$\dim_{\mathbbm{k}}\operatorname{Ext}^1(S,S)=2$
for every cyclic simple $S$.
In particular, when $m=n=1$, we have
$T=\mathbbm{k}[z]$ and
\[
 \dim_{\mathbbm{k}}\operatorname{Ext}^1
 \bigl(V_0(0,0),V_0(0,0)\bigr)=2.
\]

\begin{theorem}\label{thm:all-ext}
All first extensions between finite-dimensional simple Yetter-Drinfeld
$T$-modules are as follows: unequal $z^t$-eigenvalues give zero Ext;
truncated pairs are given by Theorem~\ref{thm:A-class} and
\eqref{eq:m1}; cyclic pairs by Proposition~\ref{prop:B-noniso}
and Theorem~\ref{thm:B-class}.
\end{theorem}

\section{Wild representation type for nonzero eigenvalues}

Fix $m>1$ and $\lambda\in\mathbbm{k}^\times$. Let
\begin{equation}\label{eq:Clambda}
 \mathcal C_\lambda=\{W\in{}_T^T\mathcal{YD}_{\mathrm{fd}}:
 (z^t-\lambda)^NW=0\text{ for some }N\ge1\}.
\end{equation}
Since $m=\prod_i e_i>1$, fix an index $i$ with $e_i>1$.
For the cyclic simple $S=V_\lambda(0,0)$,
\eqref{eq:di} gives $d_j=1$ for all $j$.

For arbitrary $N\ge1$ and matrices
$A=(a_{kj}),B=(b_{kj})\in M_N(\mathbbm{k})$, define $E(A,B)$
with generators $w_{0,j},w_{1,j}$ ($1\le j\le N$) and relations
\begin{equation}\label{eq:E-relations}
 gw_{0,j}=w_{0,j},\quad gw_{1,j}=w_{1,j},\quad
 (z^t-\lambda)w_{0,j}=0,\quad
 (z^t-\lambda)w_{1,j}=w_{0,j}.
\end{equation}
Its explicit PBW basis is
\begin{equation}\label{eq:E-basis}
 \{z^\ell y^{\mathbf r}w_{\varepsilon,j}:
 0\le\ell<t,\ 0\le r_i<e_i,\ \varepsilon=0,1,\ 1\le j\le N\}.
\end{equation}
Thus $\dim_{\mathbbm{k}}E(A,B)=2nN$.
The action of $g$ is given by \eqref{eq:g-on-basis}
with $a=b=0$, and the $y_i$-actions are given by
\eqref{eq:B-middle-action} with $u=1$, applied to each
pair $w_{0,j},w_{1,j}$.
Define the generator coactions by
\begin{equation}\label{eq:E-coaction}
\begin{aligned}
 \delta(w_{0,j})={}&1\otimes w_{0,j},\\
 \delta(w_{1,j})={}&1\otimes w_{1,j}
 +\sum_{k=1}^N a_{kj}z\otimes z^{t-1}w_{0,k}\\
 &+\sum_{k=1}^N b_{kj}\sum_{s=1}^{e_i-1}
 \frac{C^{(i)}_{e_i-1,s-1}(0,0)}
 {[s]_{\xi^{-tm_i^2/m_0}}}
 y_i^s\otimes z^{t-1}y_i^{e_i-s}w_{0,k}.
\end{aligned}
\end{equation}
The coaction of every other PBW vector is obtained from
\eqref{eq:YD}; this is well defined by the next lemma.

\begin{lemma}\label{lem:E-YD}
For all $A,B\in M_N(\mathbbm{k})$, formulas
\eqref{eq:E-relations}--\eqref{eq:E-coaction} define a
$2nN$-dimensional Yetter-Drinfeld module, and
$(z^t-\lambda)^2E(A,B)=0$.
\end{lemma}

\begin{proof}
Put $D=z^t-\lambda$. As a left $T$-module,
\[
 E(A,B)\cong
 \bigoplus_{j=1}^N M(0,0)/D^2M(0,0),
 \qquad w_{0,j}=Dw_{1,j}.
\]
By \eqref{eq:PBWM}, the vectors in \eqref{eq:E-basis}
form a basis. Hence
\[
 \dim_{\mathbbm{k}}E(A,B)=2nN,
 \qquad D^2E(A,B)=0.
\]

Let $L=\bigoplus_{k=1}^N Tw_{0,k}\cong S^{\oplus N}$.
Its coaction is independent of $A,B$, and
\[
 \delta(z^\ell w_{0,k})
 =1\otimes z^\ell w_{0,k}.
\]
Since $\Delta(z)=1\otimes z+z\otimes1$,
the $z$-terms in \eqref{eq:E-coaction}
satisfy coassociativity.

For the $y_i$-terms, fix $i$ and put
$q=\xi^{-tm_i^2/m_0}$. The mixed coefficients
reduce to
\begin{equation}\label{eq:E-Gaussian}
\begin{aligned}
&\frac{C^{(i)}_{e_i-1,s+r-1}(0,0)}
{[s+r]_q}
{s+r\brack s}_q\\
&\qquad=
\frac{C^{(i)}_{e_i-1,s-1}(0,0)}
{[s]_q}
C^{(i)}_{e_i-s,r}(0,0),
\end{aligned}
\end{equation}
for $s,r\ge1$ and $s+r<e_i$.
This follows from \eqref{eq:C} and the Gaussian
factorial identities, as in
Proposition~\ref{prop:B-normal}.
At $s+r=e_i$, the boundary equation holds since
\[
 C^{(i)}_{e_i-s,e_i-s}(0,0)=0
 \qquad(d_i=1).
\]
Thus all scalar coassociativity equations hold.
Since $\delta(L)$ is independent of $A,B$,
these equations are linear in $a_{kj},b_{kj}$;
no products of matrix entries occur.
Consequently, they hold for arbitrary $A,B$.
The counit identity is immediate.

For the $g$-compatibility, the additional
terms have respective weights
\[
 \xi^{-te_1m_1/m_0}=1,\qquad
 \xi^{-((t-1)e_1m_1+e_im_i)/m_0}=1,
\]
by the root compatibility. Hence
\[
 g\triangleright\delta(w_{\varepsilon,j})
 =\delta(w_{\varepsilon,j})
 \qquad(\varepsilon=0,1).
\]

Finally, under the adjoint--diagonal action,
\[
 z\triangleright(c\otimes v)
 =c\otimes zv+(zc-cz)\otimes v.
\]
All first tensor factors in \eqref{eq:E-coaction}
commute with $z$, so
\[
 \delta(z^tw_{1,j})
 =(\operatorname{id}\otimes z^t)\delta(w_{1,j}).
\]
Since $D$ annihilates $L$, we obtain
\[
 \delta(Dw_{1,j})=1\otimes w_{0,j},
 \qquad
 \delta(D^2w_{1,j})=0.
\]
Together with the $g$-compatibility, these
identities show that the defining relations
are preserved. Thus the generator coactions
extend to a well-defined coaction satisfying
\eqref{eq:YD}, and the coassociativity and
counit identities extend by the Hopf identities.
\end{proof}

\begin{proposition}\label{prop:E-hom}
For $A,B\in M_N(\mathbbm{k})$ and $A',B'\in M_{N'}(\mathbbm{k})$, every
Yetter-Drinfeld morphism
$f:E(A,B)\to E(A',B')$ is uniquely specified by
$P,R\in M_{N'\times N}(\mathbbm{k})$ satisfying
\begin{equation}\label{eq:E-intertwine}
 A'P=PA,\qquad B'P=PB,
\end{equation}
via
\begin{equation}\label{eq:E-morphism}
\begin{aligned}
 f(w_{0,j})&=\sum_{k=1}^{N'}p_{kj}w'_{0,k},\\
 f(w_{1,j})&=\sum_{k=1}^{N'}p_{kj}w'_{1,k}
 +\sum_{k=1}^{N'}r_{kj}w'_{0,k}.
\end{aligned}
\end{equation}
Consequently, $E(A,B)\cong E(A',B')$ precisely when $N=N'$
and the two matrix pairs are simultaneously similar.
\end{proposition}

\begin{proof}
Put $D=z^t-\lambda$ and use primes for the corresponding
objects in the target. By \eqref{eq:E-relations},
\[
 \ker D=\operatorname{Im}D
 =\bigoplus_{j=1}^N Tw_{0,j}\cong S^{\oplus N}.
\]
Since $fD=D'f$, the lower family is preserved.
By the $g$-action on the PBW basis and
Lemma~\ref{lem:PBW-weights}, the $g$-eigenspace of
eigenvalue $1$ in $E(A',B')$ is spanned by
$w'_{0,k}$ and $w'_{1,k}$, $1\le k\le N'$. Hence there exist
matrices $P,Q,R$ such that
\[
\begin{aligned}
 f(w_{0,j})&=\sum_k p_{kj}w'_{0,k},\\
 f(w_{1,j})&=\sum_k
 (q_{kj}w'_{1,k}+r_{kj}w'_{0,k}).
\end{aligned}
\]
The identity $D'f(w_{1,j})=f(w_{0,j})$ gives $Q=P$.

Comparing the coefficients of
$z\otimes z^{t-1}w'_{0,k}$ and
$y_i\otimes z^{t-1}y_i^{e_i-1}w'_{0,k}$
in
\[
 \delta f(w_{1,j})
 =(\operatorname{id}\otimes f)\delta(w_{1,j})
\]
yields $A'P=PA$ and $B'P=PB$, respectively.
Here $e_i>1$ ensures that $z$ and $y_i$ are
distinct PBW elements, and the coefficient
of the $s=1$ term in \eqref{eq:E-coaction}
is $1$.

Conversely, these intertwining equations
make \eqref{eq:E-morphism} action- and
coaction-preserving for arbitrary $R$.
On each pair of PBW-coordinate blocks,
$f$ has matrix
\[
 \begin{pmatrix}P&R\\0&P\end{pmatrix}.
\]
Thus $f$ is invertible precisely when
$N=N'$ and $P$ is invertible, proving
the simultaneous-similarity criterion.
\end{proof}

\begin{theorem}\label{thm:wild}
If $m>1$, then for every $\lambda\in\mathbbm{k}^\times$ the category
$\mathcal C_\lambda$ is wild in the representation-embedding sense.
More precisely, there exists a $\mathbbm{k}$-linear, exact, faithful functor
\[
 F:\operatorname{mod}_{\mathrm{fd}}\mathbbm{k}\langle X,Y\rangle
 \longrightarrow\mathcal C_\lambda
\]
which preserves indecomposability and reflects isomorphisms.
Its image consists of modules annihilated by $(z^t-\lambda)^2$.
\end{theorem}

\begin{proof}
A finite-dimensional $\mathbbm{k}\langle X,Y\rangle$-module
is specified by a pair of matrices $(A,B)$, and its
morphisms are the matrices $P$ satisfying
$A'P=PA$ and $B'P=PB$.
Define $F(A,B)=E(A,B)$ and let $F(P)$ be
\eqref{eq:E-morphism} with $R=0$, setting
$F(0)=0$.

By Proposition~\ref{prop:E-hom}, this defines a
$\mathbbm{k}$-linear functor which is faithful
and reflects isomorphisms. In the PBW basis
\eqref{eq:E-basis}, $F(P)$ acts as $P$ on each
of the $2n$ coordinate blocks, so $F$ is exact.
Moreover, Lemma~\ref{lem:E-YD} gives
$(z^t-\lambda)^2F(A,B)=0$, hence the image
of $F$ lies in $\mathcal C_\lambda$.

To prove preservation of indecomposability,
suppose $(A,B)$ is indecomposable. By
Proposition~\ref{prop:E-hom}, an idempotent
endomorphism of $E(A,B)$ corresponds to
$(P,R)$, where $P$ commutes with $A,B$ and
\[
 P^2=P,\qquad RP+PR=R.
\]
Indecomposability of $(A,B)$ forces
$P=0$ or $P=I_N$. In either case the second
equation gives $R=0$. Thus $E(A,B)$ has
no nontrivial idempotent endomorphisms and
is indecomposable. Conversely, a common
direct-sum decomposition of $(A,B)$ induces
one of $E(A,B)$.

Therefore $F$ is a representation embedding,
and $\mathcal C_\lambda$ is wild in the
asserted sense.
\end{proof}

\begin{remark}
The functor need not be full, since the morphisms in
Proposition~\ref{prop:E-hom} admit an independent matrix
$R$. Put $D=z^t-\lambda$ and
\[
 L=\ker D=\operatorname{Im}D\cong S^{\oplus N}.
\]
Every simple submodule and simple quotient of $E(A,B)$
is annihilated by $D$, while
$E(A,B)/L\cong S^{\oplus N}$. Hence
\[
 \operatorname{soc}E(A,B)
 =\operatorname{rad}E(A,B)=L.
\]
Thus $E(A,B)$ has Loewy length $2$ and
composition length $2N$.

In fact, the image of $F$ lies in the Serre
subcategory generated by $S=V_\lambda(0,0)$.
No assertion is made for the zero-eigenvalue
component or the case $m=1$.
\end{remark}

\section*{Acknowledgements}
The author would like to express sincere gratitude to
Professor Gongxiang Liu for his valuable guidance and
helpful suggestions.
The author is supported by the National Natural Science
Foundation of China (NSFC, Grant No.~12201164).

\vspace{1em}

\noindent\textsc{School of Mathematics, Hangzhou Normal University,
Hangzhou 311121, China}

\noindent\textit{E-mail address:}
\texttt{zhengyiwei@hznu.edu.cn}

\end{document}